\documentclass[a4paper,10pt,oneside,reqno]{amsart}
\usepackage{amsmath}
\usepackage{amsthm}
\usepackage{amssymb}
\usepackage{graphicx}

\theoremstyle{plain}
\newtheorem{teore}{Theorem}[section]
\newtheorem{defin}[teore]{Definition}
\newtheorem{lem}[teore]{Lemma}
\newtheorem*{notat}{Notation}
\newtheorem{coro}[teore]{Corollary}
\newtheorem{propo}[teore]{Proposition}
\newtheorem*{claim}{Claim}
\newtheorem*{fact}{Fact}
\newtheorem*{claim*}{Claim}
\newtheorem*{teor*}{Theorem}
\newtheorem*{theorem*}{Theorem}
\theoremstyle{remark}
\newtheorem{ejemplo}[teore]{{\sc Example}}
\newtheorem{notas}[teore]{{\sc Remark}}

\newcommand{\nrm}[1]{\|#1\|}
\DeclareMathOperator{\conv}{conv}
\newcommand{\norm}[1][\cdot]{\Vert #1\Vert}

\newcommand{\prop}{\begin{propo}}
\newcommand{\fprop}{\end{propo}}
\newcommand{\cor}{\begin{coro}}
\newcommand{\fcor}{\end{coro}}

\newcommand{\defi}{\begin{defin}\rm}
\newcommand{\fdefi}{\end{defin}}
\newcommand{\eje}{\begin{ejemplo}}
\newcommand{\feje}{\end{ejemplo}}
\newcommand{\lema}{\begin{lem}}
\newcommand{\flema}{\end{lem}}
\newcommand{\teor}{\begin{teore}}
\newcommand{\fteor}{\end{teore}}
\newcommand{\teorsn}{\begin{teor*}}
\newcommand{\fteorsn}{\end{teor*}}
\newcommand{\nota}{\begin{notas}\rm}
\newcommand{\fnota}{ \end{notas}}
\newcommand{\clam}{\begin{claim}}
\newcommand{\fclam}{\end{claim}}
\newcommand{\clams}{\begin{claim*}}
\newcommand{\fclams}{\end{claim*}}

\newcommand{\lclam}{\begin{lclaim}}
\newcommand{\flclam}{\end{lclaim}}
\newcommand{\prucl}{\prue[Proof of Claim:]}
\newcommand{\fprucl}{\fprue}
\newcommand{\ben}{\begin{enumerate}}
\newcommand{\een}{\end{enumerate}}
\newcommand{\bit}{\begin{itemize}}
\newcommand{\eit}{\end{itemize}}

\newcommand{\mc}[1]{\mathcal{#1}}

\newcommand{\mr}[1]{\mathrm{#1}}

\newcommand{\casos}{\begin{itemize}}
\newcommand{\fcasos}{\end{itemize}\setcounter{cs}{1}}

\newcommand{\fami}{(\mc F_i,\theta_i)_{i=1}^{r}}
\newcommand{\famii}{(\mc F_i,\theta_i)_{i\in I}}

\newcommand{\famin}{(\mc F_i,\theta_i)_{i\in I}}

\newcommand{\famib}{(\mc B_i,\theta_i)_{i=1}^{r}}
\newcommand{\famic}{(\mc C_i,\theta_i)_{i=1}^{r}}

\newcommand{\fin}{\textsc{FIN}}
\newcommand{\pe}{\preceq}

\newcommand{\uu}{\theta}
\newcommand{\ts}{\mathcal{T}}
\newcommand{\nn}{\mathbb{N}}

\newcommand{\conj}[2]{ \{ {#1}\,:\,{#2} \} }

\newcommand{\om}{\omega}

\newcommand{\buit}{\emptyset}

\newcommand{\ga}{\gamma}

\DeclareMathOperator{\conc}{^\smallfrown}

\DeclareMathOperator{\spr}{spread}

\newcommand{\Ga}{\Gamma}

\newcommand{\al}{\alpha}
\newcommand{\be}{\beta}
\newcommand{\de}{\delta}
\newcommand{\De}{\Delta}
\newcommand{\la}{\lambda}
\newcommand{\ini}{\sqsubseteq}

\newcommand{\sig}{\sigma}
\newcommand{\vphi}{\varphi}

\newcommand{\vep}{\varepsilon}

\newcommand{\N}{{\mathbb N}}

\newcommand{\rest}{\negmedspace\negmedspace\upharpoonright\negthickspace}

\newcommand{\e}{\varepsilon}

\newcommand{\supp}{\mathrm{supp\, }}

\newcommand{\ran}{\mathrm{ran\, }}

\newcommand{\con}{\subseteq}

\newcommand{\cones}{\varsubsetneq }

\newcommand{\prue}{\begin{proof}}
\newcommand{\fprue}{\end{proof}}

\makeindex

\begin{document}
\title{A classification of Tsirelson type spaces}
\author{J. Lopez-Abad}\address{ Equipe de Logique Math\'{e}matique \\
Universit\'{e} Paris 7- Denis Diderot\\
C.N.R.S. -UMR 7056\\
2, Place Jussieu- Case 7012\\
75251 Paris Cedex 05 \\
France}\email{abad@logique.jussieu.fr}
\author{A. Manoussakis}\address{Department of Mathematics\\ University of Aegean\\
Karlovasi, Samos, GR 83200, Greece}\email{amanouss@aegean.gr}
\thanks{Research partially supported by E$\Pi$EAK program
Pythagoras}
\begin{abstract}
We give a complete classification of mixed Tsirelson spaces $T[\fami]$ for finitely many pairs of given
compact and hereditary families $\mc F_i$ of finite sets of integers and $0<\theta_i<1$ in terms of the
Cantor-Bendixson indexes of the families $\mc F_i$,  and $\theta_i$ ($1\le i\le r$). We prove that there are
unique countable ordinal $\al$ and $0<\theta<1$ such that every block sequence of $T[\fami]$ has a
subsequence equivalent to a subsequence of the natural basis of the $T(\mc S_{\om^\al},\theta)$. Finally, we
give a complete criterion of comparison in between two of these mixed Tsirelson  spaces.
\end{abstract}

\maketitle

\section*{Introduction}
The line of research  we continue in this paper has been initiated by an old problem of S. Banach asking if
every Banach space contains a subspace isomorphic to $c_0$ or some $\ell_p$. This problem was solved
negatively by B. S. Tsirelson \cite{tsi}  who provided the first example of a Banach space that does not
contain any of the spaces $c_0$, $\ell_{p}$, $1\leq p<\infty$. The idea of Tsirelson's construction became
particularly apparent after T. Figiel and W. B. Johnson \cite{fig-joh} have shown that the norm  of the dual
of Tsirelson space satisfy the following implicit equation
\begin{equation}\label{dfsffbhrtetr}
\norm[\sum_{n}a_{n}e_{n}]=\max\{\sup_{n}\vert a_{n}\vert,\frac12\sup\sum_{i=1}^{d}\nrm{E_i(\sum_{n
}a_{n}e_{n})}\},
\end{equation}
where  the sequences $(E_i)_{i=1}^d$ considered above consists on successive subsets of integers with the
property that $d\le \min E_1$, $d\in\N$, and  $E_i(\sum_{n}a_ne_n)=\sum_{n\in E_i}a_n e_n $ is the
restriction of $\sum_{n}a_ne_n$ on the set $E_i$.  We refer to \cite{cas-shu} for an extended study of
Tsirelson space $T$. A first systematic abstract study on Tsirelson construction was given by  S. Bellenot
\cite{bel} and S. A. Argyros and I. Deliyanni \cite{arg-del0}. Given a real number $0<\theta<1$ and  an
arbitrary \emph{compact} and \emph{hereditary} family $\mc F$ of finite sets of integers one  defines the
Tsirelson type Banach space $T(\mc F,\theta)$ as the completion of $c_{00}$ with the implicitly given norm
\eqref{dfsffbhrtetr} replacing $1/2$ by $\theta$ and using sequences $(E_i)_{i}$ of finite sets of integers
which are $\mc F$-admissible,
 i.e. there is some
$\{m_i\}_{i=1}^d\in \mc F$ such that $m_1\le \min E_1\le \max E_1<m_2\le \min E_2\le \max E_2<\dots < m_d\le
\min E_d\le \max E_d$. In this notation, Tsirelson example is the space $T(\mc S,1/2)$, where $\mc
S=\conj{s\con \N}{\#s\le \min s}$ is the so called \emph{Schreier family}. It  was proved in \cite{arg-del0}
that if the \emph{Cantor-Bendixson index} $\iota(\mc F)$ and $\theta$ satisfy the inequality $\theta \cdot
\iota(\mc F)>1$, then the space $T(\mc F,\theta)$ is reflexive. Moreover, in the case of $\iota(\mc F)\geq
\omega$, they proved that the space $T(\mc F,\theta)$ does not contain any of the classical spaces $c_0$ or
$\ell_{p}$, $1\leq p<\infty$. In the case that $\mc{F}$ is chosen to be the family of the finite subsets of
$\N$ with cardinality at most $n\ge 2$, denoted by $[\N]^{\leq n}$, it was shown in \cite{bel},
\cite{arg-del0} that the corresponding space $T([\N]^{\leq n}, \theta)$ is isomorphic to $c_0$ if
$n\theta\leq 1$ and is isomorphic to $\ell_p$ ($1<p<\infty$) if $\theta= n^{-1/q}$, where $q$ is the
conjugate of $p$ (i.e. $1/p+1/q=1$).

Further examples of Tsirelson type spaces with interesting properties are the spaces
$T(\mc{S}_{\al},\theta)$ considered in \cite{arg-asp},\cite{arg}, where the compact and hereditary
families $\mc S_{\al}$ are the  \emph{$\al$-Schreier families}, the natural generalizations of the
Schreier family to index $\om^{\al}$ ($\mc S_1=\mc S$). These spaces share many properties with the
original Tsirelson space, and their natural Schauder bases are examples of $w-$null sequences with
large oscillation indexes. A basic property of any $\mc S_{\al}$ is that it is spreading (see
definition below). This is used to show that  every normalized block sequence with respect to their
natural bases $(e_n)$ is equivalent to a subsequence of $(e_n)$, a property that  $c_0$ and
$\ell_p$ also have. From this, and the fact that the  Cantor-Bendixson indexes of the families $\mc
S_{\al}$ and $[\N]^{\le n}$ are very much different, it can be explained why $T(\mc S_{\al},1/2)$
does not contain isomorphic copies of $\ell_p\cong T([\N]^{\le n},n^{-1/q})$ or $c_0\cong
T([\N]^{\le n}, 1/n)$.

The aim of this paper is to understand in these terms the so called \emph{mixed Tsirelson spaces}
$T[\fami]$, whose norms are defined implicitly by
\begin{equation*}
\norm[x]_{(\mc{F}_i,\uu_i)_{i=1}^r}= \max \big\{\|x\|_\infty, \sup \conj{\uu_i\sum_{j=1}^n
\|E_jx\|_{(\mc{F}_i,\uu_i)_{i=1}^r}}{ (E_j)_{j=1}^n \text{ is $\mc F_i$-admissible, $1\le i\le r$}
} \big\},
\end{equation*}
for  arbitrary compact and hereditary families $\mc F_i$ and   establish  a criterion of
comparability in between them. The first step in this direction was done by J. Bernues and I.
Deliyanni \cite{ber-del} and J. Bernues and J. Pascual \cite{ber-pas} who proved the following two
results:

a) If the Cantor-Bendixson indexes of the families are finite then  $T[\fami]$  is saturated  by
  either to $c_0$ or to some $\ell_{p}$, $1<p<\infty$.

b) If the Cantor Bendixson index of $\mc F$ is equal to $\omega+1$ then $T(\mc F,\theta)$  contains a
subspace isomorphic to a subspace of $T(\mc S,\theta)$.

The only case  left is   when  one of the families has infinite index.  Recall that every ordinal $\al>0$ has
a unique decomposition as $\al=\om^{\be}k+\de$, where $\de<\om^{\be}$ and $k\in \N$. Using it twice it
follows that every infinite ordinal $\al$ has the unique decomposition $ \al=\om^{\om^{\ga}n+\xi}m +\eta$
(see \cite{sier}).  Now given a compact family $\mc F$, let $\ga(\mc F)$ and $n(\mc F)$ be $\om^{\om^\ga}$
and $n$ in the previous decomposition for $\al$ equal to the Cantor-Bendixson index of $\mc F$. Following
this notation, our main result is the following
\begin{theorem*}
Fix $(\mc{F}_i,\theta_i)_{i=1}^r$ such that at least one of the families has infinite index. Then
there is some $1\le i_0\le r$  such that for every compact and hereditary family $\mc G$ the
following are equivalent.

\noindent (1) $\ga(\mc G)=\ga(\mc F_{i_0})$.

\noindent (2) Every infinite dimensional closed subspace of $T[\fami]$ contains a subspace
isomorphic to a subspace of $T(\mc G,\theta_{i_0}^{n(\mc G)/n(\mc F_{i_0})})$.

 \noindent (3) Every normalized block sequence $(x_n)$ of $T[\fami]$ has a subsequence
$(x_n)_{n\in M}$ equivalent to the subsequence $(e_{\min \supp x_n})_{n\in M}$ of the basis of $T(\mc
G,\theta_{i_0}^{n(\mc G)/n(\mc F_{i_0})})$.

\end{theorem*}
It readily follows that

\noindent (c) every normalized block sequence $(x_n)$ of $T[\fami]$ has a subsequence $(x_n)_{n\in
M}$ equivalent to the subsequence $(e_{\min \supp x_n})_{n\in M}$ of the basis of $T[\fami]$.

\noindent (d) There are unique countable ordinal $\al$ and  $0<\theta<1$ such that every normalized block
sequence with respect to the basis $(e_n)$ of $T[\fami]$ has a subsequence equivalent to a subsequence of the
basis $(e_n)$ of $T(\mc S_{\om^\al},\theta)$.

So, for example $T(\mc S_{\om^{3}4+\om 5},1/2^4)$ and $T(\mc S_{\om^3},1/2)$ are mutually saturated, while
$T(\mc S_{\om^3},1/2)$ and $T(\mc S_{\om^4},1/2)$ are totally incomparable.

Another consequence is that     every subspace of   $T[\fami]$  contains a $\mc S_{\om^\al}-\ell_{1}$
spreading model, that is, there exits a constant $K>1$ such that for every sequence of coefficients
$(a_n)_{n}$
$$
\norm[\sum_{n\in s}a_{n}x_{n}]\geq\frac{1}{K}\sum_{n\in s}\vert a_{n}\vert \quad  (s\in \mc
S_{\om^\al})$$ In particular, every subspace of $T[\fami]$ contains an asymptotic
$\ell_{1}$-subspace. Asymptotic $\ell_{1}$-spaces, the structure of these spaces as well as the
structure of the spreading models of a Banach space is a current research topic, which provides
interesting examples and structural results in Banach space theory (see \cite{aost}, \cite{OT1}).

The proofs given in this paper use four main ingredients: We work with the equivalent reformulation of the
implicit norm of $T[\fami]$ given by the norming set $K(\fami)$, and the so-called tree analysis of a
functional of  $K(\mc F,\theta )$  (see section \ref{section3}). In particular, given a normalized block
sequence $(x_n)$ of the basis $(e_n)$ we provide an algorithm to estimate the norm of a linear combination
$\sum_n a_n x_n$ in terms of a corresponding linear combination of a subsequence of the basis $(e_n)$ of an
auxiliary space $T[(\mc G_i,\theta_i)_{i=1}^r]$, much in the spirit of well-known works in this field.
Secondly, we use the well know fact (see \cite{gas}),\cite{arg-tod}) that given two compact and hereditary
families $\mc F$ and $\mc G$ there is an infinite set $M$ such that either $\mc F\rest N=\conj{s\in \mc
F}{s\con N}\con \mc G\rest N=\conj{s\in \mc G}{s\con N}$ or viceversa. This is indeed a consequence of the
fact that for every compact and hereditary family $\mc F$ there is an infinite set $M$ such that $\mc F\rest
M$ is, what we call here, \emph{homogeneous} on $M$. It turns out that the $\con$-maximal elements of such
families have the \emph{Ramsey property}, which we will use here to avoid some combinatorial computations.

Finally, we reduce the study of $T[\fami]$ for compact and hereditary families $\mc F_i$'s to  the
case of $T(\mc G,\theta)$ for some \emph{regular} family $\mc G$, i.e. a compact hereditary family
$\mc G$ that is in addition \emph{spreading} (see below). This additional  regularity property of
families $\mc G$ have two main advantages; the first is that the associated norming set $K(\mc
G,\theta)$ has a simpler form; the second one is that their Cantor-Bendixson index  is preserved if
we restrict them to an infinite set.

The paper is organized as follows: In the first section we   introduce notation,  basic
combinatorial definitions, and mixed Tsirelson spaces. In the second section we study the behavior
of subsequences of the natural basis of $T[\fami]$ in the case of regular families. An important
outcome  of this section is the reduction we make from finitely many families to one.

The third section is devoted to an abstract study of compact and hereditary families of finite sets
of integers.  In particular, we introduce homogeneous and uniform families and we prove two
combinatorial results, basic tools for this work. This section provide us the link  between  mixed
Tsirelson spaces built by compact and hereditary families  with Tsirelson type spaces constructed
using a regular family.

In the last section  we show that every block sequence of a mixed Tsirelson space $T[\fami]$ has a
further subsequence equivalent to a subsequence of its basis. As a consequence of this  and of the
results of the previous sections we provide several saturation results.  We give also, using
special convex combinations, two criteria to obtain incomparability for  Tsirelson type spaces.
Finally, we expose the classification of mixed Tsirelson spaces $T[\fami]$. \medskip

 \noindent
\textbf{Acknowledgments.} We thank to S. A. Argyros and S. Todorcevic for their useful remarks.
\section{Basic facts}
 Thorough all this paper we are going to deal with families of finite sets of integers. The family of all
finite sets of integers is denoted here by $\fin$. Given $s,t\in
\fin$ we write $s<t$ (resp. $s\le t$) to denote that $\max s<\min t$
(resp. $\max s\le \min t$), and for an integer $n$ we write $n<s$
($n\le s$) whenever $\{n\}<s$ (resp. $\{n\}\le s$). These orders can
be easily extended to vectors $x,y\in c_{00}(\N)$: $x<y$ ($x\le y$)
iff $\supp x <\supp y$ (resp. $\supp x\le \supp y$), where for $x\in
c_{00}$, $\supp x=\{n\in\N: x(n)\ne 0\}.$ We say that a sequence
$(s_n)$ of finite sets of integers is a \emph{block sequence} if
$s_n<s_{n+1}$ for every $n$. In a similar manner one defines the
corresponding notion of \emph{block sequence of vectors} of
$c_{00}$.

Given an infinite set $M$ and a finite set $s$ we denote
$M/s=\conj{n\in M}{n>s}$, and for a given integer $n$, let
$I/n=I/\{n\}$. The \emph{shift} of a non-empty set $A$ of integers
is ${_*}A=A\setminus\min A$.   Given two sets $A$ and $B$ we set
$A\setminus B=\conj{n\in A}{n\notin B}$, and $M\setminus
m=\conj{n\in M}{n\ge m}$.  For a given family $\mc F\con \fin$, an
infinite set $M\con \N$ and a finite set $s$, let $\mc F\rest
M=\conj{s\in \mc F}{s\con M}$ be the \emph{restriction} of $\mc F$
in $M$, and let $\mc F_{s}=\conj{t\in \fin}{s<t,\, s\cup t\in \mc
F}$. Given a finite set $s$ we use $\#s$ to denote its cardinality.
Finally, every time we write an enumeration $A=\{m_i\}$ of a set $A$
we mean an strictly increasing enumeration.

Concerning now in topological aspects, observe that  the family of all finite sets of integers has
the natural topology induced by the product topology on the Cantor space $\{0,1\}^\N$, simply by
identifying every finite set $s$ with its characteristic function $\xi_s:\N\to \{0,1\}$. We say
then that a family $\mc F\con \fin$ is \emph{compact} if $\mc F$ is closed with respect to the
previous topology. This means that there is no infinite sequence $(s_n)\con \mc F$ such that $s_n
\cones s_{n+1}$.  Given a compact family $\mc F$, recall that $\partial\mc F$ is the set of all
proper accumulation points of $\mc F$ and that $\partial^{(\alpha)}(\mc
F)=\bigcap_{\be<\alpha}\partial(\partial^{(\be)}(\mc F))$. The rank is well defined since ${\mc F}$
is countable and therefore a scattered compactum, so the sequence $(\partial^{(\al)}(\mc F))_\al$
of iterated derivatives must vanish. We define,  as in \cite{ber-del},  the Cantor-Bendixson index
$\iota(\mc F)$ of a compact family $\mc F$ as the minimal ordinal $\al$ such that
$\partial^{(\al)}\mc F\con \{\buit\}$. Observe that this definition is a slight variation of the
standard one, where one considers the first ordinal $\al$ such that $\partial^{(\al)}\mc F$
vanishes. Let us point out   the reason to take this definition of the index  of a family $\mc F$:
while for families with infinite index the results we present in this paper have exactly the same
form using the standard notion of Cantor-Bendixson index, for families with finite index the
standard Cantor-Bendixon index  cannot be used to  characterize the corresponding mixed Tsirelson
spaces (see \cite{ber-del}).

    A family $\mc F$ is called \emph{hereditary} iff $s\con t\in \mc F$
implies that $s\in \mc F$. Another relevant order of $\fin$ is $\pe$: Given two finite sets $s$ and
$t$ we write $s \preceq t$ iff $|t|= |s|$ and the only strictly increasing map $\sigma:t
\rightarrow s$ satisfies that $n \ge \sig(n)$ for all $n\in t$, or equivalently   if
$s=\{s_1,\ldots, s_{d}\}$ then $t=\{t_1,\ldots,t_d\}$ and $s_{i}\leq t_{i}$ for every $i\leq d$.
 We say that a family $\mc F$ of finite subsets of an infinite set $M$
is \emph{spreading} on $M$ if $s \preceq t \con M $ and $s\in \mc F$
implies $t\in \mc F$.
 We say that $\mc F$ is spreading if it is
spreading on $\N$.
  We say that $\mc F$ is \emph{regular} on $M$ iff
it is compact hereditary and spreading on $M$, and that $\mc F$ is regular if it is regular on
$\N$.

Examples of regular families are  the families of subsets of $M$ with cardinality $\le n$, denoted by
$[M]^{\le n}$, and with index $n$.  Indeed, we will see that every regular family with finite index is, when
restricted to some tail $\N/n$, of this form (see Proposition \ref{reg=hom}).  A regular family of index
$\om$ is the well-known \emph{Schreier family}
 $$
\mc S=\conj{s\in \fin}{\# s\le \min s}.
 $$
In general, for a countable ordinal $\al$ we can define inductively on $\al$ an
\emph{$\al$-Schreier family} by $\mc S_1=\mc S$, $\mc S_{\al+1}=\conj{s_1\cup \dots \cup
s_n}{(s_i)\con \mc S_{\al} \text{ is $\mc S$-admissible}}$, and if $\al$ is a limit ordinal, $\mc
S_\al=\bigcup_n \mc S_{\al_n}\rest (\N\setminus n)$, where $(\al_n)$ is a fixed increasing sequence
of ordinals with limit $\al$.  It can be shown  that $\mc S_\al$ is a regular family with index
$\om^\al$ \cite{arg-asp}. We introduce now two well known  operations between families of finite
sets.

\defi
Fix two families $\mc F$ and $\mc G$  of finite sets. Recall the following from \cite{arg-del}
\begin{align*}
\mc F\oplus\mc G = & \conj{s\cup t}{s<t,\, s\in \mc G,\, t\in \mc F} ,\\
 \mc F\otimes\mc G = & \conj{s_1\cup \dots \cup s_n}{(s_i) \text{ is
block, }s_i\in \mc F\text{ and }\{\min s_i\}\in \mc G }.
\end{align*}
\fdefi
    The  operation  $\mc F\oplus \mc G$ is called \emph{block
sum} while the operation $\mc F\otimes \mc G$ is called \emph{convolution}. Observe that $\al+1$-Schreier
families are defined inductively by  the formula $\mc S_{\al+1}=\mc S_{\al}\otimes \mc S$. Also, it is well
known that the index of the families $\mc F\oplus\mc G$ and $\mc F \otimes \mc G$ are equal $\iota (\mc
F)+\iota (\mc G)$ and $\iota (\mc F)\iota (\mc G)$ respectively, assuming that $\mc F,\mc G$ are regular (see
Proposition \ref{reg=hom}). So, if $\al$ has \emph{Cantor normal form} $\al=\om^{\al_0}n_0+\dots
+\om^{\al_k}n_k$ (see \cite{sier} for standard properties of ordinal arithmetic), the regular family $(\mc
S_{\al_0}\otimes [\N]^{\le n_0})\oplus \dots \oplus (\mc S_{\al_k}\otimes [\N]^{\le n_k})$ is
 of index $\al$.

It is not difficult to prove that $\oplus$ and $\otimes$ share many properties with the addition
and multiplication of ordinals. For example, $\oplus$ and $\otimes$ are associative, and they have
the distributive law $\mc F\otimes(\mc G\oplus \mc H)=(\mc F\otimes \mc G)\oplus (\mc F\otimes \mc
H)$, while in general the two operations are not commutative or $(\mc F\oplus \mc G)\otimes \mc
H\neq (\mc F\otimes \mc H)\oplus (\mc G\otimes \mc H)$ (as for the addition and multiplication of
ordinals).

In order to keep the notation easier we introduce the following notation
\begin{notat}
By $(\mc F_i,\theta_i)_{i\in I}$ we shall  mean a sequence of pairs of compact and hereditary
families $\mc F_i$ and real numbers $0<\theta_i<1$ ($i\in I$).  We call a sequence $(\mc
F_i,\theta_i)_{i\in I}$ \emph{regular} if in addition every $\mc F_i$ is regular. Given two
sequences $(\mc F_i,\theta_i)_{i\in I}$ and $(\mc F_i,\theta_i)_{i\in J}$ we use $(\mc
F_i,\theta_i)_{i\in I}\conc (\mc F_i,\theta_i)_{i\in J}$ to denote the concatenation of the two
sequences $(\mc F_i,\theta_i)_{i\in I\cup J}$.  Given $\mc F\con \fin$ and $m\in \N$ let
$$\mc F^{\otimes(m)}=\mc F \otimes \overset{(m)}{\cdots}\otimes \mc F.$$
\end{notat}

We are now ready to give the definition of mixed Tsirelson spaces.
\defi
Given a sequence $(\mc F_i,\theta_i)_{i\in I}$  the norm $\|\cdot\|_{\famin}$ on $c_{00}$ is defined as
follows. For $x\in c_{00}$ let
\begin{equation}\label{dfdfwerjkgjkgh}
\norm[x]_{(\mc{F}_i,\uu_i)_{i\in I}}= \max \big\{\|x\|_\infty, \sup \conj{\uu_i\sum_{j=1}^n
\|E_jx\|_{(\mc{F}_i,\uu_i)_{i\in I}}}{ (E_j)_{j=1}^n \text{ is $\mc F_i$-admissible, $i\in I$} }
\big\}\,.
\end{equation}
Next,  $T[\famin]$ denotes the completion of $(c_{00},\|\cdot\|_{(\mc{F}_i,\uu_{i})_{i=1}^{r}})$.
Observe that a Tsirelson type space $T(\mc F,\theta)$ is nothing else but the mixed Tsirelson space
$T[(\mc F,\theta)]$.
\fdefi

\nota \noindent (a) From the hereditariness of the families $\mc
F_i$ ($i\in I$) it follows easily that the  Hammel standard basis $(e_n)$ of $c_{00}$ is an 1-unconditional
normalized Schauder basis of $T[\famin]$. In the sequel whenever we consider block sequences will be with
respect the basis $(e_n)_n$.

\noindent (b) The basis  $(e_n)$ is also boundedly complete, and if there exists $i\in I$ with
$\theta_i>1/\iota (\mc F_i)$ (with the convention $1/\iota (\mc F_i)=0$ for $\iota (\mc F_i)$ is infinite)
then $(e_n)$ is also shrinking.  Therefore in this case $T[\famin]$ is reflexive (see \cite{arg-tod} for more
details).


\noindent (c) Observe  that if in the previous definition of the norm $\|\cdot\|_{(\mc{F}_i,\uu_i)_{i\in I}}$
we do not impose that $\mc F_i$ are necessarily  hereditary but only $\sqsubseteq$-hereditary ($s\sqsubseteq
t$ if $s\con t$ and $s< t\setminus s$) then in the corresponding completion $T[\famin]$ the sequence $(e_n)$
is still a  bimonotone Schauder basis, not necessarily unconditional.

\noindent (d) It can be shown that the implicit formula (\ref{dfdfwerjkgjkgh}) remains true for every $x\in
T[\famin]$  (see \cite{OS} or Remark \ref{fgrelthuuliw} below).

\noindent (e) If we  allow to some of the families $\mc F_i$ to be  non-compact, i.e. some of their closures
contain an infinite set, then it follows easily that $T[(\mc F_i,\theta_i)_{i=1}^{r}]$ is $\ell_1$-saturated.
Indeed,  every seminormalized block sequence contains a further subsequence, for which every finite initial
subsequence is $\mc F_i$-admissible for a non-compact family $\mc F_i$, and hence equivalent to the natural
basis of $\ell_{1}$.
 \fnota
  We present now  an standard alternative description of the norm of
the space $T[\famin]$, closer to the spirit of Tsirelson's original definition. Let us denote by
$K((\mc F_i,\theta_i)_{i\in I})$ the minimal subset of $c_{00}$

 \noindent(a) containing $\pm e^*_n$ ($n\in\nn$)

 \noindent (b) it  is closed under the \emph{$(\mc{F}_i,\uu_i)$}-operation
 ($i\in I$): $\uu(f_1+\dots+f_n)\in K(\famii)$ for every $\mc F_i$-admissible sequence $ (f_i)_{i=1}^n \con
K(\famii)$.

The norm induced by $K(\famin)$, i.e.
$$ \|x\|_{K(\famin)}=\sup \conj{f(x)}{f\in K(\famin)},\text{ for $x\in c_{00}$,}
 $$ is exactly the norm $\|x\|_{\famin}$ defined above. Given an infinite set $M$ of integers we
 set
 $$K^M(\famin)=\conj{\phi\in K(\famin)}{\supp \phi\con M\}}.$$

\nota\label{fgrelthuuliw}

\noindent (a) It is easy to see that the closure under the pointwise convergence topology of ${\conv
K(\famin)}$ is the unit dual ball $B_{T[\famin]^*}$. It  follows that $B_{T[\famin]^*}$ is closed under the
$(\mc F_i,\theta_i)$-operation ($i\in I$).

\noindent (b) For every infinite set $M$ of integers and every sequence $(a_n)_{n\in M}$ of scalars
we have
$$\nrm{\sum_{n\in M}a_n e_n}_{\famin}=\nrm{\sum_{n\in M}a_ne_n}_{K^M(\famin)}.$$
Observe that $K^M((\mc F_i,\theta_i)_{i\in I})=K((\mc F_i\rest M,\theta_i)_{i\in I})$ if $\mc F_i$
is regular for every $i\in I$, but that in general the previous inequality is not true.
 \fnota
Notice that, by minimality of $K(\famin)$, every functional from $ K(\famin)$ either has    the  form $\pm
e_n^*$ ($n\in \N$), or it is the result of a $(\mc F_i,\theta_i)$-operation to some sequence in $K(\famin)$
and $i\in I$. This suggests that somehow every element of $K(\famin)$ has a complexity, that increases in
every use of  the $(\mc F_i,\theta_i)$-operations. This is captured by the following definition.
\defi \cite{arg-del}  A family
$(f_t)_{t\in \mc T}\con K(\famin)$ is called a \emph{tree analysis} of a functional $f\in K((\mc
F_i,\theta_i)_{i\in I})$  if the following are satisfied:

\noindent (i) $\ts=(\ts,\pe_\ts)$  is a finite tree with  a unique root denoted by $\emptyset$, and
$f_\emptyset=f$.

\noindent (ii) For every $t\in \ts$ maximal node, $f_t=\e_t e^*_{k_t}$ where $\e_t=\pm 1$.

\noindent (iii) For every $t\in\ts$ which is not maximal, there exists  $i\in I$ such that
$(f_s)_{s\in S_t}$ is $\mc{F}_i$-admissible and $f_t=\uu_i \sum_{s\in S_t} f_s$, where $S_t$
denotes  the set of immediate successor nodes of $t$.

Note that $S_t$   is well ordered by $s_0<s_1$ iff $\supp f_{s_0}<\supp f_{s_1}$. Whenever there is
no possible confusion we will write $\pe$ in order to denote $\pe_\mc T$.

\fdefi It is not difficult to see, by the minimality of the set  $K(\famin)$,
that every  functional of $ K(\famin)$ admits a tree analysis.

As we mentioned before in Remark  \ref{fgrelthuuliw},  in general it is not true that $K^M((\mc
F_i,\theta_i)_{i\in I})=K((\mc F_i\rest M,\theta_i)_{i\in I})$ for a given infinite set $M$ of integers,  so,
a priori,  it does not suffice to control the restrictions $\mc F_i\rest M$ ($i\in I$) for the understanding
of norms $\nrm{\sum_{n\in M}a_n e_n}_{\famin}$.  We will see soon that the following is a key definition for
this purpose.

\defi
 Given a family $\mc F$ we define the family of all $\mc F$-admissible sets as follows:
We say that a finite set $t=\{m_i\}_{i=0}^{k-1}$ \emph{interpolates} the block sequence
$(s_i)_{i=0}^{k-1}$ of finite sets iff

\begin{equation*}
m_0\le s_0<m_1\le s_1<\dots <m_{n-1}\le s_{n-1}.
\end{equation*}
We say that $t=\{n_i\}$ \emph{interpolates} $s=\{m_i\}$ iff  $t$
interpolates the block sequence $(\{m_i\})$.

Given a family $\mc F$ of finite sets,  a block sequence
$(s_i)_{i=0}^{n-1}$ of finite sets  is \emph{$\mc F$-admissible} if
there is some $t\in \mc F$ which interpolates $(s_i)_{i=0}^{n-1}$.
We define
\begin{equation*}
\mathrm{Ad}(\mc F)=\conj{\{m_i\}_{i=0}^n\in \fin}{ (\{m_i\})_{i=0}^n
\text{ is $\mc F$-admissible}},
\end{equation*}
the family of all $\mc F$-admissible finite sets. \fdefi Notice that
if $M\subseteq \N$ and $(s_i)$  is an $\mc F$-admissible sequence of
subsets of $M$, then $\{\min s_i\}\in \mathrm{Ad}(\mc F)\rest M$.
The converse is not true in general.

 We list some
properties of the $\mc F$-admissible. Particularly interesting is
the characterization of spreadness of a family in terms of its $\mc
F$-admissible sets.

\prop\label{sprewewww} \noindent (a) $\mc F\con
\mathrm{Ad}(\mc F)$.

\noindent (b) If $\mc F$ is compact or hereditary,  then so is
$\mathrm{Ad}(\mc F)$.

\noindent (c)   $\mc F$ is spreading on $M$ iff $\mathrm{Ad}(\mc
F)\rest M=\mc F$.

\noindent (d)  Set $\mathrm{Ad}^{(n+1)}(\mc F)=\mathrm{Ad}(\mathrm{Ad}^{(n)}\mc F)$, $\mathrm{Ad}^{(0)}(\mc
F)=\mc F$. Then $\spr({\mc F})=\conj{s}{  \exists t\in \mc F \,( t\pe s)}=\bigcup_n \mathrm{Ad}^{(n)}(\mc F)$
is the minimal spreading family on $\N$ containing $\mc F$. In case that $\mc F$ is compact or hereditary
then so is $\spr (\mc F)$, and if $\mc F$ is regular on some set $M$, $\spr (\mc F)\rest M=\mc F$.
 \fprop
 \prue (a), (b)  are easily proved.

(c): If $\mc F$ is spreading on $M$, and $t\in \mc F$ interpolates
some $s\con M$, then, in particular, $t\pe s$ and hence $s\in \mc
F$. Suppose that $\mathrm{Ad}(\mc F)=\mc F$, and suppose that $s\pe
t$, with $s\in \mc F$ and $t\con M$. Set $s=\{m_i\}_{i=1}^k$ and
$t=\{n_i\}_{i=1}^k$.   For each $0\le j\le k$ let
$t_j=\conj{m_i}{1\le i \le k-j}\cup \conj{n_i}{k-j+1\le i\le k}$.
Observe that $t_0=s\in \mc F$, $t_j$ interpolates $t_{j+1}$ and that
$t_k=t$, so an easy inductive argument  finishes the proof of (c).
(d) follows by similar arguments  than $(c)$.
 \fprue

Finally, let us recall the following from \cite{gas}
\teor \label{tgas1}  Suppose that $\mc F$  and $\mc G$
are two compact and hereditary families. Then there is some infinite set $M$ such that either $\mc
F\rest M\con \mc G\rest M$ or $\mc G\rest M\con \mc F\rest M$.
\fteor
As for regular families $\mc F$ we have that $\iota(\mc F\rest M)=\iota(\mc F)$ for every $M$ (see
Proposition \ref{reg=hom}), it follows that if $\mc F$ and $\mc G$ are two regular families with
$\iota(\mc F)<\iota(\mc G)$ then for every $M$ there is some $N\con M$ such that $\mc F\rest N \con
\mc G\rest N$. In other words, strict inequalities between indexes of regular families imply,
modulo restrictions, strict inclusion between those families.

\section{Subsequences of the basis for regular families.}
The purpose of this section is to understand, for regular families, the relationship between the operations
$\oplus$ and $\otimes$ on regular families and corresponding norming sets. For example, what is the relation
between $K(\mc F\oplus \mc F,\theta)$ and $K(\mc F,\theta)$? It is well known that if the family $\mc F$ has
finite index, then these two norming sets are, in general, different, as the corresponding Tsirelson type
spaces are isomorphic to different $\ell_p$'s. However if $\mc F$ is, for example, the Schreier family $\mc
S$ then it  can be easily shown that $[\N]^{\le 3}\otimes \mc S\con \mc S\otimes [\N]^{\le 2}$, and hence
$$[\N]^{\le 8}\otimes ( \mc S \otimes [\N]^{\le 2})\con ([\N]^{\le 3}\otimes ([\N]^{\le 3} \otimes \mc S))\otimes [\N]^{\le 2}     \con
(\mc S \otimes [\N]^{\le 4})\otimes [\N]^{\le 2}= \mc S\otimes [\N]^{\le 8}.$$ It follows, by
induction on the complexity of $\phi\in K(\mc S\otimes [\N]^{\le 2},\theta)$ that
$\phi=\phi_1+\dots+\phi_8$ for some block sequence $(\phi_i)_{i=1}^8\con K(\mc S, \theta)$. This
clearly implies that
$$\nrm{\sum_n a_n e_n}_{(\mc S,\theta)}\le \nrm{\sum_n a_n e_n}_{(\mc S \otimes [\N]^{\le 2},\theta)}\le 8 \nrm{\sum_n a_n e_n}_{(\mc S,\theta)}$$
for every $0<\theta<1$ and every sequence $(a_n)$ of scalars.   As one can guess this reasoning cannot be
applied  to an arbitrary regular family $\mc F$ with infinite index since we do not have an explicit
presentation of $\mc F$ as for the Schreier family. However, we do have the index of the family, and by the
properties of the ordinals we have that
$$\iota([\N]^{\le 3}\otimes (\mc F\otimes [\N]^{\le 2}))=3(\iota(\mc F) 2)< \iota(\mc F)2+\om \le \iota(\mc F)3$$
and, since $\mc F$ is regular, there is some infinite set $M$ of integers such that $[M]^{\le 3}\otimes (\mc
F\rest M \otimes [M]^{\le 2})\con \mc F\otimes [\N]^{\le 3}$, hence
$$\nrm{\sum_{n\in M} a_n e_n}_{(\mc F,\theta)}\le \nrm{\sum_{n\in M} a_n e_n}_{(\mc F \otimes [\N]^{\le 2},\theta)}\le 3 \nrm{\sum_{n\in M} a_n e_n}_{(\mc F,\theta)},$$
so the two subsequences $(e_n)_{n\in M}\con T(\mc F,\theta)$ and $(e_n)_{n\in M}\con T(\mc F\otimes [\N]^{\le
2},\theta)$ of the corresponding natural bases are  3-equivalent.

We start with the following simple fact that readily follows from the definitions of the norms.
\begin{fact}\label{lemma3.10}
Suppose that  $(\mc F_i,\theta_i)_{i\in I}$,  $(\mc G_i,\theta_i)_{i\in I}$ and $M\con \N$  have the property
 that
\begin{center}
every $\mc G_i$-admissible sequence of subsets of $M$ is $\mc F_i$-admissible ($i\in I$).
\end{center}
Then for every sequence $(a_n)_{n\in M}$ of scalars
$$
\norm[\sum_{n\in M}a_ne_n]_{(\mc G_i,\uu_i)_{i\in I}} \leq \norm[\sum_{n\in M}a_ne_n]_{(\mc
F_i,\uu_i)_{i\in I}}.$$
 \end{fact}
The next  is a simple generalization of the above fact that   will be used repeatedly.
\prop\label{xxsllki}   Suppose that  $(\mc F_i,\theta_i)_{i\in I}$,  $(\mc G_i,\theta_i)_{i\in I}$, $M\con \N$ and $k\in \N$  have
the property that
\begin{equation}\label{vfdgrnsww}
[M]^{\le k}\otimes \mathrm{Ad}(\mc F_i)\rest M\con \mc G_i\rest M \otimes [M]^{\le k} \quad (i\in
I).
\end{equation}
Then for every   sequence $(a_n)_{n\in M}$ of scalars
\begin{equation*}
\nrm{\sum_{n\in M}a_n e_n}_{\famin}\le k \nrm{\sum_{n\in M}a_n e_n}_{(\mc G_i,\theta_i)_{i\in I}}.
\end{equation*}
\fprop
 \prue
We are going to show, using (\ref{vfdgrnsww}),  that for every $\phi\in K^M(\famin)$ there are
$\psi_0<\cdots <\psi_{l-1}$ in $K^M((\mc G_i,\theta_i)_{i\in I})$, $l\le k$, such that
$\phi=\psi_0+\dots +\psi_{l-1}$. The proof is by induction on the \emph{complexity of $\phi$}: If
$\phi=e_n^*$, there is nothing to prove. Suppose that $\phi=\theta_i(\phi_0+\dots +\phi_n)$, where
$(\phi_i)_{i=0}^n\con K^M(\famin)$ is $\mc F_i$-admissible. By inductive hypothesis find for every
$j$ a set $u_j$ of cardinality at most $k$ and a block sequence $(\psi_s)_{s\in u_j}\con K^M((\mc
G_i, \theta_i)_{i\in I})$ such that $\phi_j=\sum_{s\in u_j} \psi_{s}$ ($j=0,\dots,n$). Observe that
since $( \phi_j)_{j=0}^n$ is $\mc F_i$-admissible, $\{\min \phi_j\}_{j=0}^n\in \mathrm{Ad}(\mc
F_i)$. Hence by our hypothesis (\ref{vfdgrnsww})
\begin{equation*}
t=\bigcup_{j=0}^n\conj{\min \psi_s}{s\in u_j}\in [M]^{\le k}\otimes (\mathrm{Ad}(\mc F_i))\rest
M\con\mc G_i\rest M \otimes [M]^{\le k}.
\end{equation*}
So there are $t_0<\dots <t_{l-1}$ in $\mc G_i\rest M$  ($l\le k$) such that $t=t_0\cup\dots \cup
t_{l-1}$. For   $0\le m\le l-1$ set
\begin{equation*}
\psi^{(m)}=\theta_i(\sum_{\min \psi_s\in  t_m} \psi_s)\in K^M((\mc G_i,\theta_i)_{i\in I}).
\end{equation*}
Then $\phi=\psi^{(0)}+\dots+\psi^{(l-1)}$, as desired.
 \fprue

As a consequence we obtain the following two results. The fist one is the general version of the
examples considered in the introduction of this section.

 \cor\label{jfhsdrevbtghgg}
Let $\famib$ and $\famic$  be regular sequences such that $\om\le \iota (\mc C_i)\le \iota (\mc B_i)\le \iota
(\mc C_i)k $  ($1\le i\le r$) for some integer $k\ge 1$. Then for every $M$ there is some $N\con M$ such that
the subsequences $(e_n)_{n\in N}$  of the  basis of $T[(\mc B_i,\theta_i)_{i=1}^{r}]$ and $T[(\mc
C_i,\theta_i)_{i=1}^{r}]$ are $2(k+1)$-equivalent.
\fcor
 \prue
By our assumption   on the indexes of the families we obtain that
$$\iota([\N]^{\le k+1}\otimes \mc B_i)=(k+1)\iota(\mc B_i)<\iota(\mc B_i)+\om \le \iota(\mc C_i
\otimes [\N]^{\le k+1})$$
for every  $1\le i\leq r$. Hence  is some $N_0\con M$ such that
$[N_0]^{\le k+1}\otimes \mc B_i\rest N_0\con \mc C_i \otimes [N_0]^{\le k+1}$ for every $1\le i\le
r$. Proposition \ref{xxsllki} yields to
\begin{equation}\label{eq25a}
\norm[\sum_{n\in N_0}a_ne_n]_{(\mc{B}_{i},\theta_i)_{i=1}^{r}} \le  (k+1)\norm[\sum_{n\in
N_1}a_ne_n]_{(\mc{C}_{i},\theta_i)_{i=1}^{r}}
\end{equation}
By Theorem \ref{tgas1}   there exists $N\con N_0$ such that
$$
[N]^{\leq 2}\otimes\mc C_i\rest N\con \mc B_i\otimes [N]^{\leq 2}\,\,\textrm{for every}\,\, i\leq
r. $$ Proposition \ref{xxsllki} yields
$$
\norm[\sum_{n\in N}a_ne_n]_{(\mc{C}_{i},\theta_i)_{i=1}^{r}} \leq 2\norm[\sum_{n\in
N}a_ne_n]_{(\mc{B}_{i},\theta_i)_{i=1}^{r}}
$$
which completes the proof.
\fprue

The next result  says the \emph{shift} operator is, when restricted to some subsequence of the basis, always
bounded. For a given set $N$ and $n\in N$, let $n^+\in N$ be the immediate successor of $n$ in $N$, i.e.
$n^+=\min N/n$.

\cor\label{jfdfsfsdhsdrevbtghgg} Let $\famib$ be a regular sequence. Then
for every $M$ there is some $N\con M$ such that for every sequence of scalars $(a_n)_{n\in N}$,
\begin{equation*}
\nrm{\sum_{n\in N}a_n e_n}_{\famib)}
   \le
\nrm{\sum_{n\in N}a_{n} e_{n^+}}_{\famib}
  \le
2\nrm{\sum_{n\in N}a_n e_n}_{\famib}.
\end{equation*}
 \fcor
  \prue
We set $I=\{1\le i\leq r\, : \, \iota(\mc B_i) \text{ is finite}\}$ and $J$ for complement of $I$. By Theorem
\ref{tgas1} we can find $N\con M$ such that and
\begin{equation*}
[N]^{\le 2}\otimes ((\mc B_i \rest N)\oplus [N]^{\le 1})\con (\mc B_i\rest N)\otimes [N]^{\le 2} \quad (i\in
J).
\end{equation*}
Moreover, we may assume that $\mc B_i\rest N=[N]^{\le \iota(\mc B_i)}$ for every $i\in I$ (see Proposition
\ref{reg=hom}). By Proposition \ref{xxsllki} we get
   \begin{equation}\label{e3.61}
\nrm{\sum_{n\in N}a_n e_n}_{(\mc B_i,\theta_i)_{i\in I}\conc(\mc B_i \oplus [\N]^{\le 1},\theta_i)_{i\in
J}}\le 2\nrm{\sum_{n\in N}a_n e_n}_{\famib}.
  \end{equation}
Observe that for every finite set $s\con N$, setting $s^+=\conj{n^+}{n\in s}\in \mc B_i$, then for $i\in I$
it holds that $s\in \mc{F}_i$, while for $i\in J$, $s^+\pe {_*}s$, hence ${_*}s\in \mc B_i$ ($\mc B_i$ is
spreading) and so $s\in \mc B_i\oplus [\N]^{\le 1}$. This fact proves that
\begin{equation}\label{e3.62}
\nrm{\sum_{n\in N}a_{n} e_{n^+}}_{\famib} \le \nrm{\sum_{n\in N}a_n e_n}_{(\mc B_i, \theta_i)_{i\in I}\conc
(\mc B_i\oplus [\N]^{\le 1},\theta_i)_{i\in J}}.
\end{equation}
Now, using that $\mc B_i$ are spreading, by \eqref{e3.61} and \eqref{e3.62} we get,
\begin{equation*}
\nrm{\sum_{n\in N}a_n e_n}_{\famib} \le \nrm{\sum_{n\in N}a_{n} e_{n^+}}_{\famib} \le
2\nrm{\sum_{n\in N}a_n e_n}_{\famib}.
\end{equation*}
  \fprue

We examine the effect of the power operation $\mc B^{\otimes(m)}$  for regular families $\mc B$ on the
corresponding norming set. We follow some of the ideas used in the proof of the corresponding result for
Schreier families (see \cite{M},\cite{OT}).

\lema\label{dewruibndrew0}
Fix an infinite set $M$ of integers, $m\in \N$ and a  regular sequence $(\mc
B_i,\theta_i)_{i=1}^r$. Then for every  sequence $(a_n)_{n\in M}$ of scalars
\begin{equation}\label{jjiifgfgs23r4230}
\nrm{\sum_{n\in M} a_n e_n}_{(\mc B_1^{\otimes(m)},\theta_1^{m})\conc(\mc B_i ,\theta_i)_{i=2}^r}
 \le \nrm{\sum_{n\in M} a_n e_n}_{ (\mc B_i,\theta_i)_{i=1}^r  }.
\end{equation}
\flema
\prue
For simplicity, using that the families considered here
 are regular, we may assume that $M=\N$.
 Suppose that $\phi\in K((\mc B_1^{\otimes(m)},\theta_1^{m})\conc(\mc B_i
,\theta_i)_{i=2}^r)$. We will show that
\begin{equation}\label{jksdfknlweriuhdfbd}
\phi(\sum_{n}a_n e_n)\le \nrm{\sum_n a_n e_n}_{(\mc B_i,\theta_i)_{i=1}^r}.
\end{equation}
It can be easily shown by induction on $m$ that if $(s_i)_{i=1}^k$ is $\mc
 B_1^{\otimes(m)}$-admissible, then
\begin{equation}
\theta^m\sum_{i=1}^k\nrm{\sum_{n\in s_i}a_ne_n}_{(\mc B_i,\theta_i)_{i=1}^r}\le \nrm{\sum_{n\in
\bigcup_{i=1}^k s_i}a_n e_n}_{(\mc B_i,\theta_i)_{i=1}^r}.
\end{equation}
It is not difficult show  by induction on the complexity of $\phi$ that the last inequality gives
\eqref{jksdfknlweriuhdfbd}.
\fprue

\lema\label{dewruibndrew} Suppose that $M$ is an infinite set   and
that $(\mc B_i,\theta_i)_{i=1}^r$ is a regular sequence  such that
\begin{equation}\label{bfnmdrjewree} \mc B_1\rest M\otimes \mc B_i\con  \mc B_i\otimes \mc B_1
\end{equation} for every $1\le i\le r$. Then for every integer $m$,
\begin{equation}\label{jjiifgfgs23r423}
\theta_1^{m-1}\nrm{\sum_{n\in M} a_n e_n}_{ (\mc B_i,\theta_i)_{i=1}^r  }   \le \nrm{\sum_{n\in M}
a_n e_n}_{(\mc B_1^{\otimes(m)},\theta_1^{m})\conc(\mc B_i ,\theta_i)_{i=2}^r}
 \le \nrm{\sum_{n\in M} a_n e_n}_{ (\mc B_i,\theta_i)_{i=1}^r  }  .
\end{equation}

\flema
\prue
The second inequality is given by the previous Lemma \ref{dewruibndrew0}. We assume that $M= \N$.
In order to prove the first inequality of (\ref{jjiifgfgs23r423}) we are going  to show that
\begin{equation} \phi(\sum_n a_n
e_n)\le\frac{1}{\theta^{m-1}}\nrm{\sum_n a_n e_n}_{(\mc B_1^{\otimes(m)},\theta_1^{m})\conc(\mc B_i
,\theta_i)_{i=2}^r)}
\end{equation}
for every if  $\phi\in K((\mc B_i ,\theta_i)_{i=1}^r )$:  For suppose that $(\phi_t)_{t\in \mc T}$
is a tree analysis of $\phi$. For every $s \pe t $ and $1\le i \le r$ let
\begin{equation}
l_i(s,t)=\#(\conj{u}{s\pe u \precneqq t \text{ and } \phi_u=\theta_i\sum_{v\in S_{u}}\phi_v}).
\end{equation}
So we have the decomposition
\begin{equation}
\phi=\sum_{t\in \mc A}( \prod_{i=1}^r \theta_i^{n_i(t)} )(-1)^{\vep_t} e_{m_t},
\end{equation}
where $\mc A$ is the set of terminal nodes of $\mc T$, $n_i(t)=l_i(\buit,t)$, $\vep_t\in \{0,1\}$,
and $m_t$ is an integer.
\clam Suppose that  there is some  $0\le d < m$ such that $ n_1(t)
\equiv  d$  ($\mod m$) for every $t\in \mc A$. Then there are $(\psi_i)_{i=1}^l\con K((\mc
B_1^{\otimes(m)},\theta_1^{m})\conc(\mc B_i ,\theta_i)_{i=2}^r)$ such that

\noindent (a) $\phi=\theta_1^{d}(\psi_1+\dots +\psi_l$)

\noindent (b) $(\psi_i)_{i=1}^l$ is $\mc B_1^{\otimes(d)}$-admissible.
\fclam
Assuming the Claim,  for every $t\in \mc A$, let $0\le d_t<m$ be such that $ n_1(t)+d_t \equiv 0$
($\mod m$), and let
\begin{equation}
\psi=\sum_{t\in \mc A}( \theta_1^{d_t}\prod_{i=1}^r \theta_i^{n_i(t)} )(-1)^{\vep_t} e_{m_t}.
\end{equation}
By the Claim we have that $\psi\in K((\mc B_1^{\otimes(m)},\theta_1^{m})\conc(\mc B_i
,\theta_i)_{i=2}^r)$. Finally,
\begin{equation}
|\phi(\sum_n a_n e_n)|\le \frac{1}{\theta_1^{m-1}}|\psi(\sum_n a_n  e_n)|\le
\frac1{\theta^{m-1}}\nrm{\sum_n a_n e_n}_{(\mc B_1^{\otimes(m)},\theta_1^{m})\conc(\mc B_i
,\theta_i)_{i=2}^r}.
\end{equation}
which completes the proof of the Lemma. \fprue
     \prucl The proof is
by induction on the complexity of $\phi$. Suppose first that $\phi=\pm e_s$. Then $d=0$ and the
desired result is clearly true. Now suppose that $\phi=\theta_j(\phi_1+\dots +\phi_k)$. There are
two cases to consider. If $j=1$, then, by inductive hypothesis applied to each $\phi_i$ ($1\le i\le
k$),  we have that for every $1\le i\le k$,
\begin{equation}
\phi_i=\theta_1^{\bar{d}}(\psi_1^{(i)}+\dots +\psi_{s_i}^{(i)})
\end{equation}
where  $0\le \bar{d}<m$ is such that $\bar{d}\equiv d-1$ ($\mod m )$,
$(\psi_l^{(i)})_{l=1}^{s_i}\con K((\mc B_1^{\otimes(m)},\theta_1^{m})\conc(\mc B_i
,\theta_i)_{i=2}^r)$ is $\mc B^{\otimes(\bar{d})}$-admissible. It follows that
\begin{align}
\phi= & \theta_1(\phi_1+\dots +\phi_k)= \left\{
\begin{array}{ll} \theta_1^m( \sum_{i=1}^k (\sum_{l=1}^{s_i}\psi_{l}^{(i)}) ) & \text{ if $d=0$} \\
\theta_1^d( \sum_{i=1}^k (\sum_{l=1}^{s_i}\psi_{l}^{(i)}) )&  \text{ if $d>0$}
\end{array}\right.
\end{align}
Using that $(\phi_i)_{i=1}^k$ is $\mc B_1$-admissible we obtain that
\begin{equation}
\bigcup_{i=1}^k\{\min \psi_j^{(i)}\}_{j=1}^{s_i}\in
\left\{\begin{array}{ll} \mc B_1^{\otimes(m)} & \text{if $d=0$, } \\
\mc B_1^{\otimes(d)} & \text{if $d>0$.}
\end{array}\right.
\end{equation}
So if $d=0$ we obtain that $\phi\in K((\mc B_1^{\otimes(m)},\theta_1^{m})\conc(\mc B_i
,\theta_i)_{i=2}^r)$, as desired;  otherwise,   (a) and (b) in the claim are clearly true for
$\phi$.

Now suppose that $j>1$. By inductive hypothesis applied to each $\phi_i$ ($1\le i\le k$), we have
that for every $1\le i\le k$,
\begin{equation}
\phi_i=\theta_1^{d}(\psi_1^{(i)}+\dots +\psi_{s_i}^{(i)})
\end{equation}
where $(\psi_l^{(i)})_{l=1}^{k_i}\con K((\mc B_1^{\otimes(m)},\theta_1^{m})\conc(\mc B_i
,\theta_i)_{i=2}^r)$ is $\mc B_1^{\otimes({d})}$-admissible. It follows that the sequence
$(\psi_1^{(1)},\dots,\psi_{s_1}^{(1)},\dots,\psi_1^{(k)},\dots,\psi_{s_k}^{(k)})$ is $ (\mc
B_1^{\otimes (d)})\otimes \mc B_j$-admissible. Observe that (\ref{bfnmdrjewree}) and the
associative property of $\otimes$ give that
\begin{equation}
(\mc B_1^{\otimes (d)})\otimes \mc B_j = (\mc B_1\otimes \overset{(d)}{\cdots} \otimes \mc
B_1)\otimes \mc B_j \con \mc B_j \otimes (\mc B_1^{\otimes (d)}),
\end{equation}
so it follows that
$(\psi_1^{(1)},\dots,\psi_{s_1}^{(1)},\dots,\psi_1^{(k)},\dots,\psi_{s_k}^{(k)})$ is also $ \mc
B_j\otimes (\mc B_1^{\otimes (d)})$-admissible.  Let $(t_i)_{i=1}^h$ be a block sequence of finite
sets such that
\begin{equation}
\conj{\min \psi_{p}^{(i)}}{1\le i\le k,\, 1\le p\le s_i}= \bigcup_{i=1}^h t_i
\end{equation}
with $\{t_i\}_{i=1}^h\con \mc B_j$ and $\{\min t_i\}_{i=1}^h\in \mc B_1^{\otimes (d)}$. For every
$1\le l\le h$ let
\begin{equation}
 \xi_l=\theta_j \sum_{ \min \psi_{p}^{(i)} \in t_l} \psi_{p}^{(i)}\in K((\mc B_1^{\otimes(m)},\theta_1^{m})\conc(\mc B_i
,\theta_i)_{i=2}^r).
\end{equation}
Whence we obtain the decomposition
\begin{equation}
\phi=\theta_1^d \sum_{l=1}^h \xi_l,
\end{equation}
giving the desired result.
\end{proof}
As a consequence of the previous lemma we get the next proposition which is  the   natural generalization of
a well know fact for the Schreier families $\mc S_n$ ($n\in \N$).
 \prop\label{powerisnot}
Let $\mc B$ be a regular family. Then for every $0<\theta<1$, every $m$, and every sequence of
scalars $(a_n)$
\begin{equation}\label{jjiifgfgs}
\nrm{\sum_n a_n e_n}_{(\mc B^{\otimes(m)},\theta^m)}
 \le
\nrm{\sum_n a_n e_n}_{ (\mc B,\theta) }
 \le
\frac{1}{\theta^{m-1}}\nrm{\sum_n a_n e_n}_{(\mc B^{\otimes(m)},\theta^m)}.
\end{equation} \qed
 \fprop

The next lemma intends to analyze the case  of indexes  $\iota (\mc B)=\om ^{\al+\be}$ and $\iota (\mc C)=\om
^{\al}$ with $\al\geq \om$ and $\be<\al$, for example $\mc C=\mc S_{\om^2+\om}$ and $\mc G=\mc S_{\om^2}$.

 \lema\label{surprise1}
 Let $M$ be an infinite set of integers, $\mc C$, $\mc B_i$ be regular families ($1\le i\le r$)
 such that $[M]^{\le 2}\con \mc C$ and
\begin{equation}\label{duirtnjkfhdww}
[M]^{\le 2}\otimes\mc C\rest M\otimes \mc B_i\rest M\con \mc B_i \otimes [\N]^{\le 2} \quad (1\le i \le r).
\end{equation}
Then for every sequence $(\theta_i)_{i=1}^{r}\subset (0,1)$ and every sequence of scalars
$(a_n)_{n\in M}$,
\begin{equation*}
\nrm{\sum_{n\in M} a_n e_n}_{(\mc B_i,\theta_i)_{i=1}^r}\le  \nrm{\sum_{n\in M} a_n e_n}_{(\mc
B_1\otimes \mc C,\theta_1)\conc (\mc B_i,\theta_i)_{i=2}^r} \le \frac{2}{\theta_1} \nrm{\sum_{n\in
M} a_n e_n}_{(\mc B_i,\theta_i)_{i=1}^r}.
\end{equation*}
\flema
\prue
The first inequality is clear. Let us show the second one. In order to keep the notation simpler,
we may assume, since all families here are regular, that $M=\N$.
\clam Every $\phi\in  K((\mc B_1\otimes \mc C,\theta_1)\conc (\mc B_i,\theta_i)_{i=2}^r  )$ has a decomposition
\begin{equation*}
\phi=\phi_1+\dots+ \phi_n,
\end{equation*}
where $(\phi_i)_{i=1}^n\con K((\mc C\otimes \mc B_1,\theta_1)\conc (\mc B_i,\theta_i)_{i=2}^r)$ is
$\mc C$-admissible.
\fclam
 \prucl
Fix $\phi\in  K((\mc B_1\otimes \mc C,\theta_1)\conc (\mc B_i,\theta_i)_{i=2}^r  )$. If $\phi=\pm
e_n^*$, the claim is clear. Now there are two cases to consider.

\noindent \textsc{Case 1.} $\phi=\theta_1(\phi_1+\dots+\phi_n)$, where $(\phi_i)_{i=1}^n\con K((\mc
B_1\otimes \mc C,\theta_1)\conc (\mc B_i,\theta_i)_{i=2}^r )$ is $\mc B_1\otimes \mc C$-admissible.
By inductive hypothesis, for each $i=1,\dots,n$,
\begin{equation*}
\phi_i=\sum_{j=1}^{n_i}\psi_j^{(i)}
\end{equation*}
where $(\psi_j^{(i)})_{j=1}^{n_i}\con K((\mc C\otimes \mc B_1,\theta_1)\conc  (\mc
B_i,\theta_i)_{i=2}^r)$ is $\mc C$-admissible, i.e. $s_i=\{\min \psi_j^{(i)}\}_{j=1}^{n_i}\in \mc
C$.
 Since for every $i=1,\dots,n$,
$\min s_i=\min \supp \phi_i$ we obtain, by (\ref{duirtnjkfhdww}), that
\begin{equation*}
s_1\cup \dots \cup s_n\in \mc C\otimes (\mc B_1\otimes \mc C)=(\mc C\otimes \mc B_1)\otimes \mc C.
\end{equation*}
Hence we can find a block sequence $(t_i)_{i=1}^m$  such that
\begin{equation*}
s_1\cup\dots \cup s_n=t_1\cup \dots \cup t_m
\end{equation*}
and such that $(t_i)_{i=1}^m\con \mc C\otimes \mc B_1$ is $\mc C$-admissible. For every $k\in
t_1\cup\dots\cup t_m$, let $i(k),j(k)$ be such that
\begin{equation*}
\min \psi_{j(k)}^{(i(k))}=k
\end{equation*}
For every $i=1,\dots,m$, let
\begin{equation*}
\psi_i=\theta_1(\sum_{k\in t_i}\psi_{j(k)}^{i(k)}).
\end{equation*}
Since  $(\psi_{j(k)}^{(i(k))})_{k\in t_i}$ is a block sequence, and since  $\{\min
\psi_{j(k)}^{(i(k))}\}_{k\in t_i}=t_i\in \mc C\otimes\mc B_1 $ we  obtain that $\psi_i\in K((\mc
C\otimes \mc B_1,\theta_1)\conc (\mc B_i,\theta_i)_{i=2}^r)$. It is clear that
\begin{equation*}
\phi=\theta_1(\phi_1+\dots+\phi_n)=
 \theta_1(\sum_{i=1}^n\sum_{j=1}^{n_i}\psi_j^{(i)})=\theta_1\sum_{i=1}^{m}\sum_{k\in
t_i} \psi_{j(k)}^{(i(k))}
 = \sum_{i=1}^{m}\theta_1\sum_{k\in t_i}
\psi_{j(k)}^{(i(k))}=\psi_1+\dots+\psi_m.
\end{equation*}
Note that $\min \psi_i=\min t_i$ ($1\le i \le m$), hence $\{\min\psi_i\}_{i=1}^m=\{\min
t_i\}_{i=1}^m\in \mc C$, so we are done.

\noindent \textsc{Case 2.} $\phi=\theta_j(\phi_1+\dots+\phi_n)$, where $(\phi_i)_{i=1}^n\con K((\mc
B_1\otimes \mc C,\theta_1)\conc (\mc B_i,\theta_i)_{i=2}^r )$ is $\mc B_j$-admissible for some
$2\le j\le r$. By inductive hypothesis, for each $i=1,\dots,n$,
\begin{equation*}
\phi_i=\sum_{j=1}^{n_i}\psi_j^{(i)}
\end{equation*}
where $(\psi_j^{(i)})_{j=1}^{n_i}\con K((\mc C\otimes \mc B_1,\theta_1)\conc  (\mc
B_i,\theta_i)_{i=2}^r))$ is $\mc C$-admissible, i.e. $s_i=\{\min \psi_j^{(i)}\}_{j=1}^{n_i}\in \mc
C$. It follows, by (\ref{duirtnjkfhdww}) and the fact that $[\N]^{\le 2}\con \mc C$,  that
\begin{equation*}
s_1\cup \dots \cup s_n\in \mc C\otimes \mc B_j \con \mc B_j\otimes [\N]^{\le 2}\con \mc B_j \otimes
\mc C.
\end{equation*}
Following similar ideas than in the proof of the Case 1 one can easily find the desired
decomposition of $\phi$.
 \fprucl
From the claim we obtain that $\theta_1\phi\in K((\mc C\otimes\mc B_1)\conc(\mc
B_i,\theta_i)_{i=2}^r)$  for every $\phi\in K((\mc B_1 \otimes \mc C,\theta_1)\conc (\mc
B_i,\theta_i)_{i=2}^r)$. Now this fact implies that for every sequence $(a_n)_n$ of scalars
\begin{equation}\label{euifgnberw1}
\nrm{\sum_{n}a_n e_n}_{(\mc B_1\otimes \mc C,\theta_1)\conc (\mc B_i,\theta_i)_{i=2}^r}\le
\frac{1}{\theta_1}\nrm{\sum_{n}a_n e_n}_{(\mc C\otimes \mc B_1,\theta_1)\conc (\mc
B_i,\theta_i)_{i=2}^r}.
\end{equation}
Since (\ref{duirtnjkfhdww}) holds, we can apply  Proposition \ref{xxsllki} to get that
\begin{equation}\label{euifgnberw2}
\nrm{\sum_{n}a_n e_n}_{(\mc C\otimes\mc B_1,\theta_1)\conc (\mc B_i,\theta_i)_{i=2}^r}\le
2\nrm{\sum_{n}a_n e_n}_{(\mc B_i,\theta_i)_{i=1}^r}.
\end{equation}
Finally we obtain the desired inequality by joining (\ref{euifgnberw1}) and (\ref{euifgnberw2}).
\fprue

\teor\label{surprise}
Suppose that  $\mc B_0$ and $\mc B_1$ are two regular families such that $\iota (\mc
B_0)=\om^{\al+\be}$, $\iota (\mc B_1)=\om^{\al}$, with $0<\be<\al$, and $\al\ge \om$. Then for
every infinite set $M$ of integers there  is an infinite  $N\subseteq M$ such that $(e_n)_{n\in
N}\con T(\mc B_0,\theta) $  and $(e_n)_{n\in N}\con T(\mc B_1,\theta)$ are equivalents.
 \fteor
\prue
Let $\mc C$ be a regular family with $\iota (\mc C)=\om^{\be}$. Since $\iota (\mc B_1\otimes \mc
C)=\omega^{\al+\be}=\iota (\mc B_0)$ passing to a subset $N$ of $M$, we may assume that the
subsequence $(e_n)_{n \in N}$ is equivalent in the  spaces $T(\mc B_0, \theta_0)$  and $T(\mc
B_1\otimes \mc C, \theta_0)$, and hence we may assume that $\mc B_0=\mc B_1\otimes \mc C$. The
result follows from the previous lemma.
\fprue

\subsection{ Reduction from finite to one}
The aim of this subsection is to reduce finite regular sequences to one, more precisely, we show in Theorem
\ref{vmnjdgdjgjrre} that for every finite regular sequence $(\mc F_i,\theta_i)_{i=1}^r$ there is some $1\le
i_0\le r$ and some infinite set $M$ of integers such that $(e_n)_{n\in M}\con T[(\mc F_i,\theta_i)_{i=1}^r]$
and $(e_n)_{\in M}\con T(\mc F_{i_0},\theta_{i_0})$ are equivalent, where $i_0$ will come from a certain
ordering of the pairs $(\mc F_i,\theta_i)$.

\defi\label{fgjrioeeee}
Recall that every ordinal $\al>0$  has the unique
decomposition
\begin{equation*}
\al=\om^ {\la(\al)} k(\al)+\xi(\al)
\end{equation*}
with  $k(\al)$   an integer and $\xi(\al)<\om^ {\la(\al)}$. Define
\begin{align*}
\ga(\al)= & \left\{ \begin{array}{ll} k(\al) & \text{if $\al$ is finite} \\
\om^{\om^{\la(\la(\al))}} & \text{if $\al$ is infinite},
\end{array} \right. \\
n(\al)= & \left\{ \begin{array}{ll} 1 & \text{if $\al$ is finite} \\
 k(\la(\al)) & \text{if $\al$ is infinite}.
\end{array} \right.
\end{align*}
For example, $\ga(\om^{\om^23+\om}4+\om^5)=\om^{\om^2}$, $n(\om^{\om^23+\om}4+\om^5)=3$ and $\ga(m)=m$ for
every integer $m$. In general for an arbitrary ordinal $\al$ we have the decomposition
\begin{equation*}
\al=\ga(\al)^{n(\al)}\om^{\xi(\la(\al)) }k(\al)+\xi(\al),
\end{equation*}
with the convention of $\xi(0)=0$.
\fdefi
We want to compare two Tsirelson type spaces $T(\mc F_0,\theta_0)$ and $T(\mc F_1,\theta_1)$. There
is  the following natural relation of domination:  we  write $(\mc F_0,\theta_0) \le '(\mc
F_1,\theta_1)$ iff there is some $C\ge 1$ such that  every subsequence $(e_n)_{n\in M}$ of the
basis of $T(\mc F_0,\theta_0)$ has a further subsequence $(e_n)_{n\in N}$  such that
$$\nrm{\sum_{n\in N}a_n e_n}_{(\mc F_0,\theta_0)}\le C\nrm{\sum_{n\in N}a_n e_n}_{(\mc F_1,\theta_1)}$$
It is clear that if $\mc F_0\con \mc F_1$ and $\theta_0\le \theta_1$ the pair $(\mc F_1,\theta_1)$
dominates $(\mc F_0,\theta_0)$. As we have already seen in Proposition \ref{powerisnot} the pairs
$(\mc F, \theta)$ and $(\mc F^{\otimes (n)},\theta_n)$ are mutually dominated ($n\in \N$). This
suggests the following more appropriate   relation:  $(\mc F_0,\theta_0)\le '' (\mc F_1,\theta_1)$
iff there are $n_0,n_1\in \N$ such that for every $M$ there is $N\con M$ such that $\mc
F_0^{\otimes (n_0)}\rest N\con \mc F_1^{\otimes (n_1)}$ and $\theta_0^{n_0}\le \theta_1^{n_1}$.

We have also shown that $(\mc S_{\om^{\al}+\be}\otimes [N]^{\le k},\theta)$, $(\mc
S_{\om^{\al}+\be},\theta)$ and $(\mc S_{\om^{\al}},\theta)$ are all of them mutually dominated,
that leads to the following definition:

\defi\label{jrjekjegd}
For pairs $(\al,\theta)$ of ordinals and real numbers we write $(\al_0,\theta_0) \le_\mathrm{T}
(\al_1,\theta_1)$   iff

\noindent (1)   $\al_0 \al_1 <\om$ and $\log_{\ga(\al_0)}\theta_0\le \log_{\ga(\al_1) }\theta_1$,
or

\noindent (2) $\al_0 \al_1 \ge \om$ and there are integers  $m_0,m_1 $ such that
$\ga(\al_0)^{n(\al_0)m_0}\le \ga (\al_1)^{n(\al_0) m_1}$ and $\theta_0^{m_0} \le  \theta_1^{m_1}$.

We write $(\mc F_0,\theta_0)\le_\mr T (\mc F_1,\theta_1)$ iff  $(\iota(\mc F_0),\theta_0)\le_\mr T (\iota(\mc
F_1),\theta_1)$.
\fdefi
To keep the notation easier  we will write $\ga(\mc F)$ for $\ga(\iota(\mc F))$ and $n(\mc F)$ for
$n(\iota(\mc F))$.  Few more properties:
\prop\label{ueersys}
\noindent (a) Suppose that $\max\{\al_0,\al_1\}\ge \om$. Then  $(\al_0,\theta_0)\le_\mr{T} (\al_1,\theta_1)$
iff $\ga(\al_0)<\ga(\al_1)$, or if $\ga(\al_0)=\ga(\al_1) $ then $\theta_0^{n(\al_1)}\le
\theta_1^{n(\al_0)}$.

\noindent (b) $<_\mathrm{T}$ is a total ordering.

\fprop \prue 
(b):  We show that $<_\mathrm{T}$ is total. So, fix two pairs $(\al_i,\theta_i)$ ($i=0,1$). Suppose first
that $\al_i \om \le \al_j $ for $i\neq j$. Then let $n$ be such that $\theta_i^n< \theta_j $. Then clearly
$\al_i n<\al_j$, and $ \theta_i^n<\theta_j $, so $(\al_i,\theta_i)<_\mathrm{T}(\al_j,\theta_j)$. Suppose now
that $\ga(\al_0)=\ga(\al_1)$. Then if $\theta_0^{n(\al_1)}\le \theta_1^{n(\al_0)}$ we obtain that
$(\al_0,\theta_0)\le_\mr T (\al_1,\theta_1)$, and  $(\al_1,\theta_1)\le_\mr T (\al_0,\theta_0)$ otherwise.

\fprue

\lema\label{nsdkfjwes} Suppose that $\Ga$ is a finite set of  countable ordinals and
$n\in \N$. There is a sequence $(\mc B_\ga)_{\ga\in \Ga}$ of regular families  such that:

\noindent (a) $\iota (\mc B_\ga)=\ga$ for every $\ga\in \Ga$.

\noindent (b) $\mc B_{\ga}=[\N]^{\le \ga}$ if $\ga\in \Ga$ is finite.

\noindent (c) For every $m_1,m_2\le n$ and every $f_i:\{1,\dots,m_i\}\to \Ga$ ($i=1,2$),
\begin{equation}
\text{if  $\prod_{i\le m_1} f_1(i)< \prod_{i\le m_2} f_2(i)$, then $\mc B_{f_1(1)}\otimes
\dots\otimes \mc B_{f_1(m_1)}\con \mc B_{f_2(1)}\otimes \dots\otimes \mc B_{f_2(m_2)}$.}
\end{equation}
\flema
 \prue Fix for every $\ga\in \Ga$  a regular family $\mc C_\ga$ of index $\ga$, with the extra requirement that if $\ga$ is finite then
 $\mc C_\ga=[\N]^{\le
\ga}$. Since $\{\mc C_{f(1)}\otimes \dots \otimes \mc C_{f(m)}\,:\, f:\{1,\dots,m\}\to \Ga,\, m\le
n\}$ is a finite set of regular families, we can find  an infinite set $M$ such that for every
$m_1,m_2\le n$ and every $f_i:\{1,\dots,m_i\}\to \Ga$ ($i=1,2$), if  $\prod_{i\le m_1} f_1(i)<
\prod_{i\le m_2} f_2(i)$, then $\mc C_{f_1(1)}\rest M\otimes \dots\otimes \mc C_{f_1(m_1)}\rest
M\con \mc C_{f_2(1)}\rest M\otimes \dots\otimes \mc C_{f_2(m_2)}\rest M$. Let $\Theta:M\to \N$ be
the unique order-preserving onto mapping between $M$ and $\N$. Then $(\Theta" (\mc C_{\ga}\rest
M))_{\ga\in \Ga}$ is the desired sequence.
\fprue

\teor \label{vmnjdgdjgjrre} Suppose that  $(\mc B_i,\theta_i)_{i=1}^r$ is a regular sequence with at least one of the
families  with infinite index. Let $i_0$ be such that $(\iota (\mc
B_{i_0}),\theta_{i_0})=\max_{<_\mathrm{T}}\conj{(\iota (\mc B_i),\theta_i)}{1\le i\le r} $. Then
every subsequence $(e_n)_{n\in M}$ of the natural basis of $T[(\mc B_i,\theta_i)_{i=1}^r)]$ has a
further subsequence $(e_n)_{n\in N}$ equivalent to the corresponding subsequence $(e_n)_{n\in N}$
of the natural basis of $T(\mc B_{i_0},\theta_{i_0})$.
\fteor
\prue To simplify the notation we assume that $M=\N$.  We re-order $(\mc B_i,\theta_i)_{i=1}^r$ in such a way that
$(\mc B_i,\theta_i)\le_\mathrm{T} (\mc B_j,\theta_j)$ for every $1\le i \le j\le r$.

Recall the decomposition (see definition \ref{fgjrioeeee})
\begin{equation}
\iota(\mc B_i)= \ga_i^{n_i}\de_i k_i+\xi_i,
\end{equation}
where if $\iota(\mc B_i)$ is finite  then  $\ga_i=\ga(\iota(\mc B_i))$,   $\de_i=n_i=k_i=1$ and
$\xi_i=0$, and if $\iota(\mc B_i)$ is infinite  then
  $n_i=n(\iota(\mc B_i))$, $\de_i=\om^{\xi(\la(\iota(\mc B_i)))}$,
$k_i=k(\iota(\mc B_i))$ and $\xi=\xi(\iota(\mc B_i))$. Observe that $\ga_r=\max\conj{\ga_i}{1\le
i\le r}$ is infinite.  Define $m_i\in \N$ ($1\le i\le r-1$) as
\begin{equation}
\label{jnkeriuhwww} m_i=\left\{\begin{array}{ll}[\log_{\theta_i}\theta_{r}]+1 & \text{if
 $\ga_i<\ga_r$} \\
n_r & \text{if $\ga_i=\ga_r$}
\end{array}\right.
\end{equation}
where $[a]$ stays for the entire part of $a$. Use the previous Lemma \ref{nsdkfjwes}  for
$\Ga=\conj{\ga_i,\de_i}{1\le i\le r} \cup \{2\}$ and $n$ large enough (for example $n= 2\max\conj{n_i m_i
}{1\le i\le r } +2$) to find the corresponding sequence $(\mc H_\ga)_{\ga\in \Ga}$ of regular families.

For $1\le i \le r$, let
\begin{align*}
\mc C_i =   (\mc H_{\gamma_i})^{\otimes (n_i)}\otimes \mc H_{\de_i}  .
\end{align*}
Observe that $\iota(\mc C_i)=\ga_i^{n_i} \om^{\de_i}$ for every $1\le i\le r$. It readily follows
that there is $N\con M$ such that for every $1\le i\le r$, if $\iota(\mc B_i)$ is infinite then
\begin{align}\label{rsweoihijhr}
[N]^{\le 2}\otimes \mc C_i\rest N\con & \mc B_i\otimes [N]^{\le 2}, \text{ and }\\
[N]^{\le k_i+1}\otimes \mc B_i\rest N\con &  \mc C_i \rest N \otimes [N]^{\le k_i+1} \nonumber
\end{align}
while
\begin{equation*}
\mc B_i\rest N=\mc C_i\rest N
\end{equation*}
 if $\iota(\mc B_i)$ is finite. Since the families $\mc B_i$ and $\mc C_i$ are regular ($1\le i \le r$),
 Proposition \ref{xxsllki}  gives that
for every sequence of scalars $(a_n)_{n\in N}$ we have that
\begin{equation}\label{kjnferijhws0}
\frac12\nrm{\sum_{n\in N}a_n e_n}_{(\mc C_i,\theta_i)_{i=1}^r}\le \nrm{\sum_{n\in N}a_n e_n}_{(\mc
B_i,\theta_i)_{i=1}^r} \le (1+ \max_{1\le i \le r,\, \iota(\mc B_i)\text{
infinite}}k_i)\nrm{\sum_{n\in N}a_n e_n}_{(\mc C_i,\theta_i)_{i=1}^r}
\end{equation}
Let $\{\varrho_i\}_{i=1}^s$ be the strictly increasing enumeration of the set $\{\ga_i\, : \, 1\le
i\le r, \, \ga_i \text{ infinite}\}$. Define
\begin{align*}
I_0=& \conj{1\le i\le r}{\ga_i \text{ is finite}} \\
I_i= & \conj{1\le j\le r-1}{\ga_j=\varrho_i} \quad (1\le i \le s),
\end{align*}
and $I_{s+1}=\{r\}$.

Finally, set $J_i= I_i\cup \dots \cup I_{s+1}$ ($0\le i\le s+1$).  The next result is the reduction
from $(\mc C_i,\theta_i)_{i=1}^r$ to $(\mc C_r,\theta_r)$.
\clam
For every $0\le j\le s$ and every sequence of scalars $(a_n)$ we have that
\begin{equation}\label{dfewtthgirtee1}
\nrm{\sum_n a_n e_n}_{(\mc C_i,\theta_i)_{i\in J_j}} \le \prod_{i\in I_j }\frac1{\theta_i^{m_i-1}}
\prod_{ i\in I_j , \,\de_i>1} \frac{2}{\theta_i} \nrm{\sum_n a_n e_n}_{(\mc C_i,\theta_i)_{i\in
 J_{j+1}}}.
\end{equation}
\fclam
\prucl
Fix $0\le j\le s$.  Let $K_j=\conj{i\in I_j}{\de_i>1}$, and suppose it is non-empty. This implies,
in particular, that $j>0$. Notice that $\varrho_j=\min\conj{\ga_k}{k\in J_j}$. So it follows that
$\de_k<\ga_k=\varrho_j\le \ga_i$ for $k\in K_j$ and $i\in J_j$. So,
\begin{align*}
2 \de_k \ga_i^{n_i}\de_i = & \ga_i^{n_i}\de_i< \ga_i^{n_i}\de_i 2 \quad (k\in K_j, \, i\in J_j) \\
2 \de_k \ga_i^{n_i}  = & \ga_i^{n_i} < \ga_i^{n_i}  2 \quad (i,k\in K_j).
\end{align*}
Hence,
\begin{align*}
[\N]^{\le 2}\otimes \mc H_{\de_k}\otimes \mc C_i \con  & \mc C_i \otimes [\N]^{\le 2}\quad (k\in
K_j,\, i\in J_j)\\
[\N]^{\le 2}\otimes \mc H_{\de_k}\otimes \mc H_{\ga_i}^{\otimes(n_i)} \con  & \mc
H_{\ga_i}^{\otimes(n_i)} \otimes [\N]^{\le 2}\quad (i,k\in K_j).
\end{align*}
A repeated application of Lemma \ref{surprise1} gives that
\begin{equation}\label{kjnferijhws1}
\nrm{\sum_n a_n e_n}_{(\mc C_i,\theta_i)_{i\in J_j}}\le \prod_{i\in K_j}\frac2{\theta_i}\nrm{\sum_n
a_n e_n}_{(\mc H_{\ga_i}^{\otimes (n_i)},\theta_i)_{i\in K_j}\conc (\mc C_i,\theta_i)_{i\in
J_j\setminus K_j}}.
\end{equation}
Using that
\begin{equation*}
\ga_k^{n_k} \ga_i^{n_i}\de_i =
  \ga_i^{n_i}\de_i<\ga_i^{n_i}\de_i\ga_k^{n_k} \quad (k\in I_j,
i\in J_{j+1})
\end{equation*}
it follows that
\begin{align*}
\mc H_{\ga_k}^{\otimes (n_k)}\otimes \mc C_i \con \mc C_i \otimes    \mc H_{\ga_k}^{\otimes (n_k)}
\quad (k\in I_j,\, i \in J_{j+1})
\end{align*}
Since it is trivial that $\mc H_{\varrho_j}^{\otimes (n_k)}\otimes \mc H_{\varrho_j}^{\otimes
(n_i)}=\mc H_{\varrho_j}^{\otimes (n_k+n_i)}=\mc H_{\varrho_j}^{\otimes (n_i)}\otimes \mc
H_{\varrho_j}^{\otimes (n_k)}$ ($i,k \in I_j$), the assumptions of   Lemma \ref{dewruibndrew} are
fulfilled, therefore
\begin{align}\label{kjnferijhws2}
 \nrm{\sum_n a_n e_n}&_{(\mc H_{\ga_i}^{\otimes (n_i)},\theta_i)_{i\in K_j}\conc  (\mc C_i,\theta_i)_{i\in J_j\setminus K_j}}= \nonumber \\
& = \nrm{\sum_n a_n e_n}_{(\mc H_{\varrho_j}^{\otimes
(n_i)},\theta_i)_{i\in I_j}\conc (\mc C_i,\theta_i)_{i\in J_{j+1}}}\le \nonumber  \\
& \le   \prod_{i\in I_j}\frac1{\theta_i^{m_i-1}}\nrm{\sum_n a_n e_n}_{(\mc H_{\varrho_j}^{\otimes
(n_i m_i)},\theta_i^{m_i})_{i\in I_j}\conc (\mc C_i,\theta_i)_{i\in J_{j+1}}}.
\end{align}
It is not difficult to see, by the choice of $m_i$'s, that the  relations
\begin{align*}
\text{$\mc H_{\varrho_j}^{\otimes (n_i m_i)}\con \mc C_r$ while $\theta^{m_i}\le \theta_r$ ($i\in
I_j$) if $j<s$ or}\\
\text{$\mc H_{\varrho_s}^{\otimes (n_i m_i)}= \mc H_{\varrho_s}^{\otimes (n_i n_r )}\con \mc
C_r^{\otimes (n_i)}$ and $\theta_i^{m_i}=\theta_i^{n_r} \le \theta_r^{n_i}$ ($i\in I_j$) if $j=s$}
\end{align*}
are true. Hence, by Lemma \ref{dewruibndrew0} in the case of $j=s$, we obtain that
\begin{equation}\label{kjnferijhws3}
\nrm{\sum_n a_n e_n}_{(\mc H_{\varrho_j}^{\otimes (n_i m_i)},\theta_i^{m_i})_{i\in I_j}\conc (\mc
C_i,\theta_i)_{i\in J_{j+1}}}\le \nrm{\sum_n a_n e_n}_{ (\mc C_i,\theta_i)_{i\in J_{j+1}}}.
\end{equation}
It is clear now that \eqref{dfewtthgirtee1} follows from equations \eqref{kjnferijhws1},
\eqref{kjnferijhws2} and \eqref{kjnferijhws3}.
\fprucl
A repeated application of previous claim gives that
\begin{equation}\label{dfgasfwe4ww}
\nrm{\sum_{n}a_n e_n}_{(\mc C_i,\theta_i)_{i=1}^r}\le \prod_{i=1}^{r-1}\frac1{\theta_i^{m_i-1}}
\prod_{ i=1 , \,\de_i>1}^{r-1}\frac{2}{\theta_i} \nrm{\sum_{n}a_n e_n}_{(\mc C_r,\theta_r)}.
\end{equation}
 It follows from  \eqref{rsweoihijhr}, \eqref{kjnferijhws0} and
\eqref{dfgasfwe4ww} that
\begin{equation*}
 \nrm{\sum_{n\in N}a_n e_n}_{(\mc B_i,\theta_i)_{i=1}^r} \le 2(1+ \max_{1\le i \le r,\, \iota(\mc B_i)\text{
infinite}}k_i)\prod_{i=1}^{r-1}\frac1{\theta_i^{m_i-1}} \prod_{ i=1 ,
\,\de_i>1}^{r-1}\frac{2}{\theta_i} \nrm{\sum_{n\in N}a_n e_n}_{(\mc B_r,\theta_r)}
\end{equation*}
\fprue
In  Theorem \ref{vmnjdgdjgjrre} we made the assumption that at least one family $\mc B_i$  has
infinite index ($1\le i\le r$). The conclusion of this theorem is also true for families all of
them with finite indexes but its proof uses different methods (see \cite{ber-del}).

\section{Topological and combinatorial aspects of families of finite sets of integers}
The main result of this section is that for every compact and hereditary family $\mc F$ there is a
regular family $\mc B$ with the same index than $\mc F$ and an infinite set $M$ of integers such
that every $\mc B$-admissible sequence of subsets of $M$ is also $\mc F$-admissible. The main tool
we use is the notion of \emph{homogeneous} family.

We start with the following list of useful properties. We leave their proofs to the reader.
 \prop\label{easy}
Fix a compact family $\mc F$, and a countable ordinal $\al$. Then

\noindent (1) For every $m\in \N$, $(\partial^{(\al)}\mc F)\rest \N/m=\partial^{(\al)}(\mc F\rest
\N/m)$.

 \noindent (2) $\emptyset\neq s\in \partial^{(\al)}\mc
F$ iff ${_*}s\in \partial^{(\al)}(\mc F_{\{\min s\}})$.

\noindent (3) For every  $n\in \N$,  $\partial^{(\al)}(\mc F_{\{n\}})=(\partial^{(\al)}\mc
F)_{\{n\}}$.

\noindent (4) $\iota (\mc F)=\al$ iff (4.1) $\partial^{(\al)}(\mc F_{\{n\}})=\buit$ for every $n$,
and (4.2) for every $\be<\al$ there is some $n$ such that $\partial^{(\be)}\mc{F}_{\{n\}}\neq
\buit$.

\noindent (5) $\iota (\mc F)=\al$ limit implies that $\partial^{(\al)}\mc F=\{\buit\}$ and for
every $n\in \N$, $\iota (\mc F\rest (\N/n))=\al$.

\noindent (6) $\iota (\mc F)=\al+1$ implies that for every  $k$ there is $m\in\N$ such that for
every $n\geq m$, $\iota (\mc F_{\{k\}}\rest \N/{n})\leq \al$.
    \qed
 \fprop

\subsection{Homogeneous families and admissible sets}
For our study it would be very useful to have a characterization of \emph{every} compact hereditary family in
terms of a class of families with good structural properties allowing inductive arguments, as for example the
Schreier families have. This is indeed the case for the class of homogeneous families. The following
definition is modelled on the  the notion of $\al$-uniform family of Pudlak and R\"{o}dl (see \cite{arg-tod}).
\defi
We say that a family $\mc F$ is \emph{$\al$-homogeneous on} $M$
($\al$ a countable ordinal) iff $\buit\in \mc F$ and,

\noindent (a) if $\alpha =0$, then $\mc F =\{\buit\}$;

\noindent (b) if  $\alpha=\beta+1$, then  $\mc F_{\{n\}}$ is
$\beta$-homogeneous on $M/n$ for every $n\in M$;

\noindent (c) if $\alpha > 0$ limit, then  there is an increasing
sequence $\{\alpha_{n}\}_{n \in M}$ of ordinals converging to
$\alpha$ such that $\mc F_{\{n\}}$ is $\alpha_{n}$-homogeneous on
$M/n$ for all $n \in M$.

$\mc F$ is called \emph{homogeneous on $M$} if it is
$\al$-homogeneous on $M$ for some countable ordinal $\al$.

Recall the following well known combinatorial notion (\cite{arg-tod}). A family $\mc F$ is
\emph{$\al$-uniform on} $M$ ($\al$ a countable ordinal) iff $\mc F=\{\emptyset\}$ for $\alpha=0$ or
$\mc F$ satisfies $(b)$ or $(c)$ in the other cases, where homogeneous is replaced by uniform. The
relationship between uniform and homogeneous families will be exposed in the Proposition
\ref{jfjgfngfnjnhjieruirt} below.
 \fdefi
\begin{notat}
If $s,t\in \fin$   we write $s\sqsubseteq t$ iff $s$ is an initial
segment of $t$.
\end{notat}
 \nota
 \noindent (a) It is easy to see that the only
$n$-homogeneous families on $M$ are the families of subsets of $M$ with cardinality $\le n$,
denoted by $[M]^{\le n}$. A well known $\om$-homogeneous family on $\N$ is the {Schreier family},
and, in general, $\om$-homogeneous families on $M$ are of the form $\conj{s\con M}{\#s\le f(\min
s)}$, with $f:M\to \N$ a unbounded and increasing mapping. Observe that all those examples are
regular families.

\noindent (b) In general,  an arbitrary homogeneous family is not
regular. However, we will show that homogeneous families are always
$\sqsubseteq$-closed, hence compact. Also, it can be shown that if
$\mc F$ is a homogeneous family on $M$, there is  $N\con M$ such
that $\mc F\rest N$ is hereditary (see \cite{arg-tod}).
 \fnota
          Uniform families  and regular families have many properties in
common. One of the most remarkable is the fact that  the index of
these families never decrease when taking restrictions. We  expose
this analogy and some others in the next proposition.
 \prop
\label{reg=hom} Suppose that $\mc F$ and $\mc G$ are homogeneous
(regular) families on $M$. Then:

\noindent(a) If $\iota (\mc F)$ is finite, then $\mc F=[M]^{\le \iota (\mc F)}$ if $\mc F$ is
homogeneous on $M$, while $\mc F\rest (M/n)=[M/n]^{\le \iota (\mc F)}$ for some $n\in M$ if $\mc F$
is regular on $M$.

\noindent(b) $\mc F_{\{n\}}$ ($n\in M$) is  homogeneous (regular) on
$M/n$.

\noindent(c) If $\mc F$ is $\al$-homogeneous, then $\partial^{(\al)}\mc F=\{\buit\}$. Hence $\iota
(\mc F)=\al$.

\noindent(d) if $\mc F$ is $\al$-homogeneous (regular) on $M$ and $N\subseteq M$ then $\mc F\rest
N$ is $\al$-homogeneous (regular) and $\iota (\mc F\rest N)=\iota (\mc F)$ for every $N\con M$.

\noindent(e) $\mc F\oplus \mc G$ and $\mc F\otimes \mc G $  are homogeneous (regular), $\iota (\mc
F\oplus \mc G)=\iota (\mc F)+\iota (\mc G)$ and $\iota (\mc F\otimes \mc G)=\iota (\mc F)\iota (\mc
G)$.

\noindent(f) If $\iota (\mc F)<\iota (\mc G)$ then for every $M$ there  is $N\con M$ such that $\mc
F\rest N \cones\mc G\rest N$.
  \fprop
\prue Suppose first that $\mc F$ is  homogeneous. (a) and (b) can be
shown by an easy inductive argument.

(c):   Now suppose first than $\al=\be+1$.  By the inductive
hypothesis, for every $n\in M$, $(\partial^{(\be)} (\mc
F))_{\{n\}}=\partial^{(\be)}(\mc F_{\{n\}})=\{\buit\}$. So,
$[M]^{\le 1}= \partial^{(\be)} (\mc F)$ (since $
\partial^{(\be)} (\mc F)$
is closed and it contains all singletons $\{n\}$ ($n\in M$)). Hence
$\partial^{(\be+1)}(\mc F)=\{\buit\}$. Suppose now that $\al$ is a
limit ordinal.  Now by the inductive hypothesis we can conclude that
for every $n\in M$,
\begin{equation}\label{njfgnuirgr}
\partial^{(\al_n)}(\mc F_{\{n\}})=(\partial^{(\al_n)}(\mc F))_{\{n\}}=\{\buit\},
\end{equation}
where $\al_n=\iota (\mc F_{\{n\}})$ is such that $(\al_n)_n$ is increasing and with limit $\al$. By
(\ref{njfgnuirgr}), $\buit \in
\partial^{(\al)}\mc F$. If there were some $s\in
\partial^{(\al)}\mc F$, $s\neq \buit$,
then $s\in \partial ^{(\al_{n}+1)}\mc F$ for every $n$, and hence  $\partial^{(\al_{n})}(\mc
F_{\{\min s\}})\ne \{\emptyset\}$, a contradiction.

(d) It follows easily by induction on $\al$  using $(a)$.

(e) is shown by induction on $\iota (\mc G)$. (f): By Proposition \ref{tgas}, there is some $N\con
M$ such that either $\mc F\rest N\con \mc G\rest N$  or else $\mc G\rest N\con \mc F\rest N$. The
second alternative is impossible since it implies that $\iota (\mc F)=\iota (\mc F\rest N)\ge \iota
(\mc G\rest N)=\iota (\mc G)$.

Finally suppose that we are dealing with regular families. (a): First, note that there must be some
$s\in \mc F$ with $|s|=\iota (\mc F)$, since otherwise, $\mc F\con [M]^{<\iota (\mc F)}$ and so,
$\iota (\mc F)<\iota (\mc F)$ an absurd. In a similar way one shows that $\mc F\con [M]^{\le \iota
(\mc F)}$. All this shows that $\mc F\rest (M/s)=[M/s]^{\iota (\mc F)}$. (b) and (c) are clear.
(d): Fix $N\con M$, and let $\Theta:M\to N$ be the unique order-preserving onto mapping between
these two sets. Since $\mc F$ is spreading on $M$, we obtain that $ \conj{\Theta"s}{s\in \mc F}\con
\mc F\rest N$; using that $\Theta "s\neq \Theta" t$ is $s\neq t$ we obtain that $\iota (\mc F\rest
N)\ge \iota (\mc F)\ge i (\mc F\rest N)$, as desired. (f)  follows from (d), while (e) is a
consequence of Theorem \ref{jfjgfewwwwpp} and (e) for homogeneous families.
 \fprue
The following result is a weaker form of  Theorem II.3.22 in
\cite{arg-tod}.
 \teor\label{jfjgfewwwwpp}
\cite{arg-tod}  Suppose that $\mc F$ is a compact and hereditary
family. Then there is some $M$ such that $\mc F\rest M$ is
homogeneous on $M$.
 \fteor
\prop\label{jfjgfngfnjnhjieruirt} Suppose that $\mc F$ is a family
of finite sets of integers. Then for every countable ordinal $\al$
the following conditions are equivalent:

\noindent (a) $\mc F$ is $\al$-homogeneous on $M$.

\noindent (b) $\mc F$ is the topological closure of an $\al$-uniform family on $M$.

\noindent (c) $\mc F$ is compact and the set $\mc
F^{\sqsubseteq-\max}$ of $\sqsubseteq$-maximal
 elements of $\mc F$ is   $\al$-uniform  on $M$. Moreover
\begin{equation}\label{jjjueweweeep}
\mc F=\conj{s\ini t}{t\in \mc F^{\sqsubseteq-\max}},
\end{equation}
hence $\mc F$  is $\ini$-hereditary, i.e. if $s\ini t\in \mc F$ then
$s\in \mc F$.
 \fprop
   \prue
(a) implies (b): The proof is   by induction on $\al$. If $\al=0$
the result is clear. Suppose that $\al=\be+1$. Then for every $n\in
M$, $\mc F_{\{n\}}$ is $\be$-homogeneous on $M/n$. Choose
$\be$-uniform families $\mc G_n$ on $M/n$ ($n\in M$) such that for
every $n\in M$, $\mc F_{\{n\}}=\overline{\mc G_n}$. Set
\begin{equation*}
\mc G=\conj{\{n\}\cup s}{s\in \mc G_n}.
\end{equation*}
It follows readily that $\mc G_{\{n\}}=\mc G_{n}$ which yields that $\mc G$ is a $\al$-uniform
family on $M$. To finish the proof we show that $\overline{\mc G}=\mc F$. First observe that if
$s\in \mc F$, $n=\min s$, then ${_*}s\in \mc F_{\{n\}}$. So, ${_*}s\in \overline{\mc G_n}$, and
hence
\begin{equation*}
s=\{n\}\cup {_*}s\in \overline{\mc G_n \oplus \{\{n\}\}}\con
\overline{\mc G}.
\end{equation*}
Now suppose that $(s_k)\con \mc G$, $s_k\to_k s\in \overline{\mc
G}$. Going to a subsequence if necessary, we may assume that $(s_k)$
is a $\De$-sequence with root $s$, i.e. $s\sqsubseteq s_k$ for every
$k$, and $(s_k\setminus s)$ is a block sequence. If $s=\buit$, then
$s \in \mc F$ by hypothesis. Otherwise, let $n=\min s$.  Then $\min
s_k=n$ for every $k$, and hence ${_*}s_k\in \mc G_{\{n\}}$. Hence
${_*}s\in \overline{\mc G_{\{n\}}}=\mc F_{\{n\}}$, and so $s\in \mc
F$. The proof if $\al$ is limit is similar.

(b) implies (c): Suppose that $\mc F=\overline{\mc G}$, where $\mc G$ is $\al$-uniform on $M$. It
is not difficult to show by induction on $\al$  that  $\mc G$ is a \emph{front on $M$} (see
\cite{arg-tod}), i.e., for every infinite $N\con M$ there is some $s\in \mc F $ such that $s
\sqsubseteq N$, and if $s, t\in \mc F$ and $s \sqsubseteq t$ then $s=t$. Observe  that the
topological closure of a front is its $\sqsubseteq$-downwards closure.
 Indeed, suppose that $s$ is a strict  initial part of some $t\in\mc
 F$. For every $m> s$ consider the set $M_{m}=s\cup M/m$.
 Using that $\mc F$ is a front
on $M$ we find $t_{m}\sqsubseteq M_{m}$ such that $t_{m}\in\mc F$,
moreover $s$ has to be initial segment of every $t_{m}$. This
implies that $t_{m}$ converges to $s$.

 So, we have that $\mc F=\conj{s\sqsubseteq t}{t\in \mc G}$. It is
clear that this implies that $\mc F^{\sqsubseteq-\max}=\mc G$.

(c) implies (a): Suppose that $\mc F$ is compact and $\mc
F^{\ini-\max}$ is $\al$-uniform on $M$. The proof is an easy
induction on $\al$ using that for every  $m\in M$, by
(\ref{jjjueweweeep}), $ \mc F_{\{m\}}=\conj{s\ini t}{t\in \mc
G_{\{m\}}}$, where $\mc G=\mc F^{\ini-\max}$.
  \fprue
The next result is the well known \emph{Ramsey} property of uniform families (see \cite{arg-tod}
for a more complete explanation of the Ramsey property).
 \prop[Ramsey Property]\label{ramseyppt}
Suppose that $\mc B$ is a $\al$-uniform family on $M$, and suppose
that $\mc B=\mc B_0\cup \mc B_1$. Then there is an infinite $N\con
M$ such that $\mc B\rest N=\mc B_0\rest N$. \qed
 \fprop
\begin{proof}
Induction on $\al$. Given $\mc B=\mc B_0\cup \mc B_1$, using
inductive hypothesis, we can find a decreasing sequence
$(M_{k})_{k}$ of infinite subsets of $M$, such that, setting
$m_{k}=\min M_{k}$, for every $k$, $M_{k+1}\con {_*}M_{k}$ and there
is an $i_{k}\in \{0,1\}$ such that $ \mc{B}_{\{m_{k}\}}\rest
M_{k+1}=(\mc{B}_{i_{k}})_{\{m_k\}}\rest M_{k+1}$. Then every $N\con
\{m_{k}\}_{k}$ for which $i_{k}$ is constant has the desired
property.
\end{proof}
 As an application of this Ramsey property   we obtain the
following two facts.
  \prop\label{tgas}
\noindent (a) \cite{arg-tod}, \cite{gas}  Suppose that $\mc F$  and $\mc G$ are two compact and
hereditary families. Then there is some infinite set $M$ such that either $\mc F\rest M\con \mc
G\rest M$ or $\mc G\rest M\con \mc F\rest M$.

\noindent (b) Suppose that $\mc F$ is homogeneous on $M$.  Then
there is some $N\con M$ such that $\mc F\rest N$ is hereditary.
\fprop
 \prue
(b):  Set $\mc B= {\mc F}^{\sqsubseteq-\max}$, and let $\mc B_0=\conj{s\in \mc B}{\mc P(s)\not\con
\mc F}$, $\mc B_1=\mc B\setminus \mc B_0$. By Ramsey, there is $N\con M$ and $i=0,1$ such that $\mc
B\rest N=\mc B_i\rest N$. If $i=1$ then we are done. Otherwise, fix $s\in \mc B\rest N$ and $t\con
s$ such that $t\notin \mc F$.  Using that $\mc F^{\sqsubseteq-\max}$ is a front on $M$, we get
$u\in \mc B\rest N$ be such that $u\sqsubseteq t \cup (N/s)$. If $t\sqsubseteq u$ then $t\in \mc
F$, which is impossible. So, $u\sqsubset t\cones s$. This means that for every $s\in \mc B\rest N$
there is some $t\cones s$, $t\in \mc B\rest N$. Hence $\buit \in \mc B\rest N$, which implies that
$\mc B\rest N=\{\buit\}$ and so $\mc F\rest N=\{\buit\}$ is hereditary.
 \fprue
   \subsection{The basic combinatorial results}
 The families $\mc F$ and $\mathrm{Ad}(\mc F)$ are in general
different, unless $\mc F$ is spreading. Nevertheless, as it is shown
in the next result, they are not so far from the topological point
of view.

\prop\label{ceyhrytyyyy} Suppose that $\mc F$ is a compact
hereditary family.   Then for every infinite set $M$ of integers
such that $\mathrm{Ad}(\mc F)\rest M$ is homogeneous on $M$,
\begin{equation*}
\iota (\mc F) \le \iota(\mathrm{Ad}(\mc F)\rest M)\le 2\iota (\mc F).
\end{equation*}
\fprop
  \prue This is done by induction on $\iota (\mc F)=\la+n$, $\la$
limit ordinal (including $\la=0$), and $n\in \N$.  Set $\mc B=\mathrm{Ad}(\mc F)\rest M$. Suppose
first that $n=0$. Observe that in this case, by Proposition \ref{easy},    for every $m\in M$ we
have that $\partial^{(\la)}(\mc F\rest \N/m)=\{\buit\}$; so for every $k$, and every $m\in M$,
$\la_k(m)=\iota (\mc F_{\{k\}}\rest (\N/m))<\la$, and $\sup \la_k(m)=\la$. Since for every $m\in M$
we have that  $\mc B_{\{m\}}=\bigcup_{k\le m}\mathrm{Ad}( \mc F_{\{k\}}\rest (\N/m) )\rest (M/m)$,
by inductive hypothesis,   we obtain that for every $m$, $\max_{k\le m} \la_k(m)\le \iota (\mc
B_{\{m\}})\le 2\max _{k\le m}\la_k(m)$. This means that $\iota (\mc B)=\la$.

Suppose now that  $\iota (\mc F)=\la+n+1$. First of all, there is some $i\in \N$  such that $\iota
(\mc F_{\{i\}})\ge \la+n$, so for every $m\in M/i$, $\iota (\mc B_{\{m\}})\ge \la+n$, and hence, by
definition of homogeneous families, $\iota (\mc B)\ge \la+n+1=\iota (\mc F)$.

Now we work to show the other inequality $\iota (\mc B)\le 2\iota (\mc F)=\la+2n+2$. We proceed by
contradiction assuming that $\iota (\mc B)\ge \la+2n+3$. By Proposition \ref{reg=hom} (d) we may
assume that $\iota (\mc B_{\{m\}})\ge \la+2n+2$ for every $m\in M$. Let $\mc G={\mc
B}^{\sqsubseteq-\max}$. Fix $m_0\in M$, and define the coloring
\begin{align*}
\Theta:\mc G_{\{m_0\}}\to \{0,\dots, m_0\}\,\,\textrm{by}&\,\,\,\
\Theta(s)=k\,\,,\,\textrm{iff}\\
&\mbox{there is some $t\in \mc F_{\{k\}}$ such that $\{k\}\cup t$
interpolates}\,\,\{m_0\}\cup s.
\end{align*} By the Ramsey property of $\mc
G_{\{m_0\}}$, we may assume, going to a subset if necessary, that $\Theta$ is constant with value
$k_0\in \{0,\dots,m_0\}$. Suppose first that  $\iota (\mc F_{\{k_0\}})\le \la+n$; then by inductive
hypothesis, $ \iota (\mc B_{\{m_0\}})\le \la + 2n$, a contradiction with our assumption. So, $\iota
(\mc F_{\{k_0\}})=\la+n+1$. Moreover, $\partial^{(\la+n+1)}(\mc F_{\{k_0\}})=\emptyset$, and hence,
$\partial^{(\la+n)}\mc F_{\{k_0\}}\not\con \{\buit\}$ is finite. Let $l=\max\partial^{(\la+n)}\mc
F_{\{k_0\}}$.
 By the Ramsey property of $\mc G_{\{m_0\}}$ we may assume that
either for every $s\in \mc G_{\{m_0\}}$ there is some $t\in \mc F_{\{k_0\}}\rest(\N/l)$ who
interpolates $s$, or else there is some $k_1\in \{m_0+1,\dots,l\}$ such that for every $s\in \mc
G_{\{m_0\}}$ there is some $t\in \mc F_{\{k_0,k_1\}}$  such that $\{k_0,k_1\}\cup t$ interpolates
$\{m_0\}\cup s$. In the first case, $\mc B_{\{m_0\}}\con \mathrm{Ad}(\mc F_{\{k_0\}}\rest (\N/l))$;
since $\iota (\mc F_{\{k_0\}}\rest (\N/l))\le \la+n$, then $\iota (\mc B_{\{m_0\}})\le \la+2n$, a
contradiction. In the second case, consider the homogeneous family $\mc B_{\{m_0,m_1\}}$, where
$m_1=\min (M/k_1)$. Then $\iota (\mc B_{\{m_0,m_1\}})\ge \la+2n+1$, and $\mc B_{\{m_0,m_1\}} \con
\mathrm{Ad}(\mc F_{\{k_0,k_1\}}\rest (\N/m_1)) $. Finally notice that $\iota (\mc
F_{\{k_0,k_1\}}\rest (\N/m_1))\le \la+n$, so $\iota (\mc B_{\{m_0,m_1\}}) \le \la +2n$, a
contradiction.
        \fprue
    \nota
The previous result is best possible. For every limit ordinal $\la$ and $n\in \N$ there is some
compact hereditary family $\mc F$ such that $\iota (\mc F)=\la+k$ and $\iota(\mathrm{Ad}(\mc
F)\rest M)=\la+2k$ for every $M$, hence $\iota(\mathrm(Ad)(\mc F)\rest M)=2\iota (\mc F)$. The
families are closely related to the example 3.10 from \cite{ber-del}. Consider a regular family
$\mc B$ on $\{2n\}$ of index $\iota (\mc B)=\la+k$, and   let $\mc F$ be the downwards closure of
$\mc F=\conj{s\cup \{k+1\}_{k\in s} }{s\in \mc B}$.  It is not difficult to prove that $\iota (\mc
F)=\iota (\mc B)=\la+k$ and that  for every $M$ $\iota((\mathrm{Ad}(\mc F))\rest M)=\la+2k$.
    \fnota
The next result somehow tell that we may assume that the given
family $\mc F$ is indeed spreading.
   \prop\label{dfderjjwewww} Fix an
arbitrary compact hereditary family $\mc F$, and an infinite set
$M$.

\noindent (a) There is some regular family $\mc B$   with the same
index than $\mc F$ and some $N\con M$ such that
\begin{equation}
\text{every $\mc B$-admissible sequence of subsets of $N$ is also
$\mc F$-admissible.}
\end{equation}
 \noindent (b) For every regular family $\mc B$ on $M$ with $\iota (\mc
B)>\iota (\mc F)$ there is some $N\con M$ such that every $\mc \overline{\mc
B}^{\sqsubseteq-\max}$-admissible sequence $(s_i)$ of subsets of $N$ with $\#s_i\ge 2$ is not $\mc
F$-admissible.
   \fprop \prue  Let $\mc C$ be a
regular and homogeneous family on $M$ such that $\iota (\mc C)\ge \iota (\mc F)=\la+r$, $\la$ limit
(including $\la=0$) and $r\in \N$.  Let $\mc G=\mc G (\mc C)=\mc C^{\sqsubseteq-\max}$. It follows,
by Proposition \ref{jfjgfngfnjnhjieruirt},  that $\mc G$ is an uniform family of $M$, as well as
$[M]^{2}\otimes \mc G$.  Since $\mc C$ is regular on $M$ with index $\ge \la+r$ we may assume that
for every $s\in \mc B$, $\# s\ge r$.  Observe that every $s\in [M]^2\otimes \mc G$ has a unique
decomposition $s=s[0]\cup \dots \cup s[l]$ with $s[0]<\dots <s[l]$, $\#s[i]=2$, and $\{\min
s[i]\}\in \mc G$, so in particular $l\ge r-1$.
 Consider the following coloring,
$$
h_{\mc C,\mc F}:[M]^2\otimes \mc G\to \{0,\dots, r, \infty\}
 $$
defined for $s\in [M]^2\otimes \mc G$ by
\begin{align*}
h_{\mc C, \mc F}(s)=k\in \{0,\dots r\}\,\,\,\,\,\,&\textrm{ iff $k$ is
minimal with the property that}\\
&\text{$(s[k],s[k+1],\dots,s[r-1],s[r+1],\dots,s[l])$ is $\mc F$-admissible},
\end{align*}
if well defined and $h_{\mc C,\mc F}(s)=\infty$ otherwise. By
Proposition \ref{ramseyppt} there is some $N=N(\mc C,\mc F)\con M$
such that $h_{\mc C,\mc F}$ is constant on $[N]^2\otimes \mc G\rest
N$ with value $k_0=k_0(\mc C,\mc F)\in \{0,\dots,r,\infty\}$.

\clam $k_0=0$ iff $\iota (\mc F)=\iota (\mc C)$. \fclam \prucl The proof is by
induction on $\iota (\mc C)$. Suppose first that $k_0=0$, and suppose that $\iota (\mc C)>\iota
(\mc F)$. By Proposition \ref{easy}, we may assume, going to an infinite subset if needed, that for
every $n\in N$, $\iota (\mc C_{\{n\}})\ge \iota (\mc F)$. Fix $n\in N$ and consider the coloring
$$
d:([N]^2\otimes (\mc G_{\{n\}})\rest N)\oplus ([N/n]^1)\to
\{0,\dots,n\}$$
 defined for $s=\{k\}\cup s[1]\cup s[2]\cup\dots\cup s[l]\in
([N]^2\otimes (\mc G_{\{n\}})\rest N)\oplus ([N/n]^1)$ by
 \begin{align*}
d(s)=j\,\,\, &\textrm{ iff there is some $t\in \mc F$ such that
$\min
t=j$ and}\\
&\textrm{ $t$ interpolates}\,\, (\{n,\min s\}, s[1],\dots, s[r-1],
s[r+1],\dots, s[l]).
\end{align*}
Observe that $d$ is well defined since we are assuming that $k_0=0$. By the Ramsey property of the
uniform family considered as domain of $d$ there is some infinite set $P\con N$ such that $d$ is
constant on $([P]^2\otimes (\mc G_{\{n\}})\rest P)\oplus ([P/n]^1)$ with value $j_0\in
\{0,\dots,n\}$. Take some $p\in P$ be such that $\iota (\mc F_{\{j_0\}}\rest (\N/p))<\iota (\mc F)$
(See Proposition \ref{easy} (6)). Then $k_0(\mc C_{\{n\}}\rest P,\mc F_{\{j_0\}}\rest (\N/p))=0$,
so, by inductive hypothesis, $\iota (\mc C_{\{n\}})=\iota (\mc F_{\{j_0\}}\rest (\N/p))) <\iota
(\mc F)$, a contradiction.

Now suppose that $\iota (\mc C)=\iota (\mc F)$.  For every $1 \le i\le r$ set $\mc C_{i}={}_{*}\mc
C_{i-1}$ and  $\mc C_0=\mc C$. Since $\mc C$ is regular  it follows easily that $\mc C_{r}$ is
regular with index $\la$. Consider the regular family $\mc D=[N]^{\le 2}\otimes  (\mc C\rest N)_r$
on $N$ with index $\iota (\mc D)=\la$. By Proposition \ref{ceyhrytyyyy}, $(\mathrm{Ad}(\mc F))\rest
N\ge \la$, so there is some $P$ such that $\mc D\rest P\con \mathrm{Ad}(\mc F)\oplus [P]^{\le 1}$.
It readily follows that ${_*}(\mc D\rest P)\con \mathrm{Ad}(\mc F)$. This shows that $k_0\in
\{0,\dots,r\}$. If $r=0$, then we are done. Now suppose that $r>0$. Let $q\in \N$ and $n_0,n_1 \in
N$ be such that $q<n_0<n_1$ and   that $\iota (\mc F_{\{q\}}\rest(\N/n_1))=\la +r-1$ (see
Proposition \ref{easy}). Since $\mc C$ is regular, we have that $\iota (\mc
C_{\{n_0\}})=\la+r-1=\iota (\mc F_{\{q\}}\rest(\N/n_1))$. By inductive hypothesis there is some
$P\con N/n_1$ such that $h_{\mc C_{\{n_0\}},\mc F_{\{q\}}\upharpoonright \,(\N/n_1)}: [P]^2\otimes
\mc G_{\{n_0\}}  \rest P\to \{0,\dots,r-1,\infty\}$ is constant with value $0$. Take arbitrary $s
\in [P]^2\otimes (\mc G_{\{n_0\}})\rest P$. Then $(s[0],\dots,s[r-2],s[r],\dots, s[l])$ is $\mc
F_{\{q\}}$-admissible, so $(\{n_0,n_1\},s[0],\dots,s[r-2],s[r],\dots, s[l])$ is $\mc F$-admissible.
Since $\{n_0,n_1\}\cup s\in [N]^2\otimes \mc G\rest N$ we obtain that $k_0=0$, as desired. \fprucl
We work now to show (a): Fix a regular and homogeneous family $\mc C$ such that $\iota (\mc
C)=\iota (\mc F)$ and such that $[\N]^{\le r}\con \mc C$.
 Let $\{n_i\}=N=N(\mc C,\mc F)$, and let
\begin{equation}
\mc D=\{s\cup t\in \mc C\rest N\, :\, s<t, \,  \# s=r \text{ and
$\exists k,l\in N$  with $s<k<l<t$ and $s\cup \{k\}\cup t\in \mc
C$}\}.
\end{equation}
It is not difficult to see that $\mc D\rest \{n_{2i}\}$  is a regular family on $\{n_{2i}\}$ with
the same index than $\iota (\mc C)= \iota (\mc F)$ (for example, by using Proposition \ref{easy}).
Now $\mc B=\spr (\mc D\rest \{n_{2i}\})$ and $\{n_{2i}\}$ have property (a): Let $(s_i)$ be a $\mc
B$-admissible sequence of finite subsets of $\{n_{2i}\}$. By Proposition \ref{sprewewww} we get
that $s=\{\min s_i\}\in \mc D\rest \{n_{2i}\}$. Write $s=t_1\cup t_2$, with $t_1<t_2$ and
$\#t_1=r$. By definition of the family $\mc D$ there are $k<l$ in $N$  such that $t_1<k<l<t_2$  and
$t_1\cup \{k\}\cup t_2\in \mc C$.  Now, since $k_0(\mc C,\mc F)=0$, it follows that $(s_i)$ is $\mc
F$-admissible.

Finally, we proceed to prove (b): Suppose that $\mc B$ is an arbitrary regular family on $M$ with
$\iota (\mc B)\ge \iota (\mc F)$. Find $M'\con M$ such that $\mc B\rest M'$ is also homogeneous on
$M'$. It is easy to see that if $(s_i)$ is a $\mc B^{\sqsubseteq-\max}$-admissible sequence of
subsets of $N(\mc B\rest M',\mc F)$,
 $\#s_i\ge 2$, which is also $\mc F$-admissible, then, since $\mc F$ is hereditary, $k_0(\mc B\rest M',\mc
F)=0$, so by the Claim, $\iota (\mc B)=\iota (\mc B\rest N')=\iota (\mc F)$. This shows (b). \fprue
  \section{block sequences of $T[(\fami)]$}\label{section3}
In this last section we show that given  finitely many compact hereditary families  $(\mc
F_i)_{i=1}^r$ such that at least one of them has infinite index there is $1\le i_0\le r$ such that
every normalized block sequence in the space $T[\famii]$ has a subsequence  equivalent to a
subsequence of the basis of the space $T(\mc F_{i_0}, \theta_{i_0})$. We first obtain this result
for the subsequences of the basis of $T[\famii]$ by applying the result of the previous section,
and in the sequel we extent this result for block sequences.

To obtain the result for a given  block sequences $(x_n)_{n}$ we show first that we can pass to a subsequence
$(x_n)_{n\in M}$   which is equivalent to the subsequence $(e_{p_n})_{n\in M}$, $p_n=\min\supp x_n$, of the
basis of the space $T([\N]^{\leq 2}\otimes \textrm{Ad}(\mc F_{i_0}), \theta_{i_0})$, for appropriate fixed
$1\le i_0\le r$. Using the results for the regular families we pass to a space $T(\mc B,\uu_{i_0})$ where
$\mc B$ is  a regular family  with $\iota (\mc B)=\iota (\mc F_{i_0})$ and moreover the subsequence
$(e_{p_n})_{n\in M}$ is equivalent in the two spaces.

Restricting the study to the families  $\mc S_{\xi}$, we obtain the  that if $(x_n)_{n}$, $(y_n)_{n}$ are
normalized block sequence in the space $T(\mc S_{\xi}, \theta)$ such that $x_n<y_n<x_{n+1}$  ($n\in \N$) then
the two sequences are equivalent.

\prop
\label{dfewekmmvv2} Fix $\fami$ with at least one of the families with infinite index. Let $i_0$ be
such that $(\mc F_{i_0},\theta_{i_0})=\max_{\le_\mr T}\{(\mc F_i,\theta_i)\}_{i=1}^r$ (See
definition \ref{jrjekjegd}). Then for every $M$ there is some $N\con M$ and a regular family $\mc
B$ with same index than $\mc F_{i_0}$ such that for every sequence $(a_n)_{n\in N}$ of scalars,
\begin{equation*}
\nrm{\sum_{n\in N}a_n e_n}_{(\mc B_{i_0}, \theta_{i_0})}\le\nrm{\sum_{n\in N}a_n e_n}_{\fami}\le
2C\nrm{\sum_{n\in N}a_n e_n}_{(\mc B_{i_0}, \theta_{i_0})}.
\end{equation*}
where the constant $C$ is given in Theorem \ref{vmnjdgdjgjrre}.
\fprop
\prue By  Proposition \ref{dfderjjwewww} we get $N_0\subseteq M$ and regular families
$\mc B_i$, with  $\iota (\mc B_i)= \iota(\mc F_i)$ ($1\le i\leq r$), such that every $\mc B_i
$-admissible sequence of subsets of $N_0$ is also $\mc F_i$-admissible.  By fact \ref{lemma3.10} it
follows that for every sequence $(a_n)_{n\in N_0}$ of scalars,
\begin{equation*}
\nrm{\sum_{n\in N_0}a_n e_n}_{(\mc B_i,\theta_i)_{i=1}^{r}}\le \nrm{\sum_{n\in N_0}a_n
e_n}_{\fami}.
\end{equation*}
Counting the corresponding indexes we can find   now  $N_1\subseteq N_{0}$ such that
 $$
[N_1]^{\leq 2}\otimes \mathrm{Ad}(\mc{F}_i)\rest N_1\subseteq (\mc {{ C}}_i\rest N_1)\otimes [N_1]^{\leq
2}\,\,\textrm{for every}\,\,i\leq r,
$$
where $\mc{{ C}}_i=\mc B_i\oplus [N_1]^{\leq 1}$ if $\iota (\mc F_i)<\omega$, $\mc{{ C}}_i=\mc B_i$
otherwise. It follows from Proposition \ref{xxsllki} that
\begin{equation*}
 \nrm{\sum_{n\in N_0}a_n
e_n}_{\fami} \le 2 \norm[\sum_{n\in N}a_ne_n]_{(\mc{{ C}}_{i},\theta_i)_{i=1}^{r}}.
\end{equation*}
 By Theorem
\ref{vmnjdgdjgjrre}, using that $\mc F_{i_0}$ has infinite index, there exist $N\subseteq N_1$ such
that
$$
\norm[\sum_{n\in N}a_ne_n]_{(\mc{{ C}}_{i},\theta_i)_{i=1}^{r}} \sim\norm[\sum_{n\in
N}a_ne_n]_{(\mc{B}_{i_0},\theta_{i_0})}.
$$
Since $\norm[\sum_{n\in N}a_ne_n]_{(\mc{B}_{i_0},\theta_{i_0})}\leq \norm[\sum_{n\in
N}a_ne_n]_{\famib}$, we get the result.
 \fprue

\nota It is worth to mention  that the conclusion of the above theorem does not hold in   case that all families $\mc F_i$'s have finite index (see
\cite{ber-del}).
\fnota

To extent the  above result to block sequences first we shall need some  preparing  work. The
following notion is  descendant of the definition of initial and partial part of a vector with
respect to a tree analysis introduced in \cite{arg-del}.

   \defi\label{jgjfjgfjfrrr}
Fix  compact and hereditary families $\mc F_i$ and real number $0<\theta_i<1$ ($i\leq r$). Let $x\in c_{00}$,
$f\in K(\fami)$ and $(f_t)_{t\in \mc T}$ a tree-analysis for $f$. Suppose that $\supp f \cap \ran x\neq
\buit$. Let $t\in \mc T$ be a $\pe$-maximal node with respect to the property that
\begin{equation*}
\supp f_t\cap \ran x=\supp f\cap \ran x.
\end{equation*}
It is clear that such $t$ exists and it is unique. Let us call it $t(x)$. Note that if $t(x)$ is not a
maximal node of $\mc T$, then, by maximality of $t(x)$, there are $s_1\neq s_2\in S_{t(x)}$ such that $\supp
f_{s_i}\cap \ran x\neq \buit$, for $i=1,2$. Observe that the set $S_t$ of immediate $\pe$-successors of $t$
is naturally ordered according to $s<t$ iff $f_s<f_t$.  Now for $t(x)$ not a maximal node, let
\begin{align*}
s_L(x)=&\min\conj{s\in S_t}{\supp f_s\cap \ran x\neq \buit},\\
s_R(x)=&\max\conj{s\in S_t}{\supp f_s\cap \ran x\neq \buit},
\end{align*}
where both minimum and maximum are with respect to the relation $<$ on $S_t$\,.

Fix now a block sequence $(x_n)_n$. For a given $n$, let $t(n)=t(x_n)$, $s_L(n)=s_L(x_n)$ and
$s_R(n)=s_R(x_n)$. For $t\in \mc T$,  we define recursively
\begin{align*}
D_t=&\bigcup_{t \preceq_\mc T u}\conj{n}{u=t(n)}=
\{n\in\N: \supp f_t\cap\ran x_n=\supp f\cap \ran x_n\}, \\
E_t=& D_t\setminus \bigcup_{s\in S_t} D_s =\conj{n}{t=t(n)}.
\end{align*}
For each $n$, set $q_n=\max\supp x_n$.
 Define recursively on $t\in
\mc T$
\begin{equation*}
g_t=\theta_i\left(\sum_{n\in E_t} \frac{f_t(x_n)}{\theta_i}e_{q_n}^*+ \sum_{s\in S_t} g_s\right)
\end{equation*}
if  $f_{t}=\uu_i\sum_{s\in S_t}f_{s}$, where $(f_{s})_{s\in S_t}$ is   $\mc{F}_{i}$-admissible.
 \fdefi
\prop\label{jfjgfruuiioo}

\noindent (a) $\supp g_t=\conj{q_n}{n\in D_t}$ for every $t\in \mc T$.

\noindent (b) The set $\{e_{q_n}^*\}_{n\in E_t}\cup \{g_s\}_{s\in S_t} $ is  $[\N]^{\le 2} \otimes
\mr{Ad}(\mc F_i) $-admissible for every  $t\in \mc T$ such that $f_{t}=\uu_i\sum_{s\in S_t}f_{s}$.
  \fprop
 \prue (a)  follows readily from the definitions.

 (b) Suppose that $\# S_t>1$, otherwise there is nothing to prove.
Let us observe that for every $s\in S_t$, $s$ not being the $<$-maximal element of $S_t$ and $D_{s}\ne
\emptyset$,
\begin{equation*}
\min\supp f_{s}\le  \supp g_{s} <\min\supp f_{s^{+}}\,.
\end{equation*}
The first inequality follows readily from (a).

Let us show now the last inequality. Assume otherwise that $\min\supp f_{s^+}\leq \max\supp
g_{s}=\max\conj{q_n}{n\in D_s} $. Then there exists $n\in D_{s}=\{m: t(m)=s\}$ such that
\begin{equation*}
\min\supp f_{s^+}\le q_n=\max\supp x_n,
\end{equation*}
hence $ \supp f_{s^+}\cap\ran x_n \neq \buit $, a contradiction since $n\in D_{s}$.  It is clear that for
every $n\in E_{t}$ such that $s_R(x_n)$ it not the $<$-maximal element of $S_t$, it holds
\begin{equation*}
\min\supp f_{s_{R}(x_n)}\le q_{n}<\min\supp f_{s_{R}(x_n)^+}.
\end{equation*}
Now setting for every $s\in S_{t}$, $B_{s}=\supp g_{s}\cup\{q_{n}: s=s_{R}(x_n), n\in E_{t}\}$, we have that
$\supp g_{t} =\cup_{s}B_{s}$ and from the previous observations we obtain that
\begin{equation*}
\min\supp f_{s}\leq   B_{s}  <\min\supp f_{s^+},
\end{equation*}
Therefore from the fact that $\{\min\supp f_{s}:s\in S_{t}\}\in \mr{Ad}(\mc{F}_i)$  we get the desired
result.
 \fprue

\prop\label{dwbvsaadfsdfsee} \noindent (a) Fix a sequence $(a_n)$ of
scalars. Then for every $t\in \mc T$,
\begin{equation*}
f_t(\sum_{n\in D_t}a_nx_n )=g_t(\sum_{n\in D_t} a_ne_{q_n}).
\end{equation*}
In particular, $f(\sum_{n}a_nx_n )=g_\buit(\sum_{n\in D_\buit} a_ne_{q_n})$.

\noindent (b) For every $t\in \mc T$, $ g_t\in \frac{1}{\theta_0}B(T[([\N]^{\le 2} \otimes \mr{Ad}(\mc
F_i),\theta_i)_{i=1}^{r}]^*)$, where $\theta_0=\min_{1\le i\leq r}\theta_i$.

\noindent (c) For every sequence $(a_n)$ of scalars
\begin{equation*}
\nrm{\sum_n a_n x_n}_{\fami}\le \frac{1}{\theta_0}\nrm{\sum_n a_n e_{q_n}}_{([\N]^{\leq 2}\otimes
\textrm{Ad}(\mc F_i),\theta_i)_{i=1}^{r}}.
\end{equation*}
 \fprop
 \prue
(a) can be shown easily by downwards induction on $t\in \mc{T}$. (b) follows from  Proposition
\ref{jfjgfruuiioo} (b) and the fact that the dual ball  of $T[([\N]^{\le 2} \otimes \mr{Ad}(\mc
F_i),\theta_i)_{i=1}^{r}]$ is closed on the $( [\N]^{\le 2} \otimes \mr{Ad}(\mc F_i),\theta_i)$-operation
(see Remark \ref{fgrelthuuliw}).

(c) Follows from (a) and (b).
  \fprue
Before we give the proof of the main result of the section we need one more auxiliary lemma.

\begin{lem}\label{b47}
Fix $(\mc F_i,\theta_i)_{i=1}^r$ with at least one of the families with infinite index, and a
normalized block sequence $(x_n)_n$   in the space $T[\fami]$. Then for every $i_0$ such that
$\iota (\mc F_{i_0})\geq \omega$ there exists infinite set $M$ such that
$$
\norm[\sum_{n\in M}a_ne_{p_n}]_{(\mc F_{i_0}, \uu_{i_0})}\leq \norm[\sum_{n\in M}a_nx_{n}]_{(\mc
F_{i}, \uu_{i})_{i=1}^r}
$$
\end{lem}

\prue
Let $(x_n)$ be a normalized block sequence    and set $p_n=\min\supp x_n$ and $P_0=\{p_n\}$. Let
$M_0$ be an infinite set of integers and let $\mc B_i$ be regular families on $\N$ with $\iota (\mc
B_i)=\iota (\mc F_i)$ such that
\begin{equation}\label{fekjddssjuuwwoo}
\text{every $\mc B_i$-admissible block sequence   of subsets of $\{p_n\}_{n\in M_0}$ is $\mc
F_i$-admissible.}\,(1\le i\le r)
\end{equation}
Let $M_1=\{m_{2i}\}$, where $\{m_i\}$ is the increasing enumeration of $M_0$.

\clam For every sequence of scalars $(a_n)_{n\in M_1}$,
\begin{align}\label{ooo1}
\nrm{\sum_{n\in M_1}a_n e_{p_n}}_{(\mc B_i,\theta_i)_{i=1}^{r}}\le & \nrm{\sum_{n\in M_1}a_n
x_n}_{(\mc F_i,\theta_i)_{i=1}^r}.
\end{align}
 \fclam
  \prucl  Choose,  for every $n$,  $\phi_n\in K(\mc F,\theta)$ such
that $\phi_n x_n \approx 1$ and $\supp \phi_n\con \supp x_n$. Define now $F:K^{M_1}((\mc B_i,
\theta_i)_{i=1}^r)\to K(\fami)$ by $F(e_{p_{n}}^*)=\phi_n$, and extend it by
$F(\theta_i(\psi_0+\dots +\psi_n))=\theta_i(F(\psi_0)+\dots F(\psi_n))$ if $(\psi_i)_{i=0}^n\con
K^{M_1}((\mc B_i, \theta_i)_{i=1}^r)$ is a $\mc B_i$-admissible block sequence ($1\le i\le r$). $F$
is well defined: Suppose that $(\psi_i)_{i=0}^n\con K^{M_1}((\mc B_i, \theta_i)_{i=1}^r)$ is  $\mc
B_i$-admissible block sequence, and set $\min \supp \psi_i=p_{m_{2k_i}}$, $\max\supp
\psi_i=p_{m_{2l_i}}$ ($ 0\le i \le n$). Then we have that for every $0\le i\le n$
\begin{equation}\label{mbkkgkttr}
\supp F(\psi_i)\con [p_{m_{2k_i}},p_{m_{2l_i+1}}]
\end{equation}
Since, by (\ref{fekjddssjuuwwoo}), $(\{p_{m_{2k_i}},p_{m_{2l_i+1}}\})_{i=0}^n$ is $\mc
F_i$-admissible the condition (\ref{mbkkgkttr}) yields that $(F(\psi_i))_{i=0}^n$ is $\mc
F_i$-admissible. It is clear now that the existence of $F$ shows the desired result.
 \fprucl
Let $i_0$ be such that $\mc F_{i_0}$ has infinite index. Then by Propositions \ref{tgas} and
\ref{ceyhrytyyyy} we can find $P\con \{p_n\}_{n\in M_1}$ such that
$$
[P]^{\leq 2}\otimes \textrm{Ad}(\mc F_{i_0})\rest P\subseteq \mc B_{i_0}\otimes [P]^{\leq 2}
$$
It follows that  $(e_{p})_{p\in P}\con T(\mc F_{i_0},\uu_{i_0})$  and  $(e_{p})_{p\in P}\con T(\mc
B_{i_0},\uu_{i_0})$ are equivalent. This, combined with the previous claim, completes the proof.
\fprue

\teor\label{infiniteindex}
Fix a finite sequence $(\mc F_i,\theta_i)_{i=1}^r$ of compact hereditary families and real numbers
such that there is some $1\le i\le r$ such that $\iota (\mc F_i)$ is infinite. Then there is $1\le
i_0\le r$  such that every normalized block sequence $(x_n)\con T[(\mc F_i,\theta_i)_{i=1}^r]$ has
a subsequence $(x_n)_{n\in N}$ which is equivalent to the subsequence $(e_{p_n})$ of the natural
basis $(e_n)$ of $T(\mc F_{i_0},\theta_{i_0})$, and where $p_n=\min \supp x_n$.

\fteor

\prue  Let $(x_n)\con T[(\mc F_i,\theta_i)_{i=1}^r]$ be a
normalized block sequence. By Proposition \ref{dwbvsaadfsdfsee} we get
\begin{equation*}
\nrm{\sum_n a_n x_n}_{(\mc F_i,\theta_i)_{i=1}^r}\le C \nrm{\sum_n a_n e_{q_n}}_{([\N]^{\leq 2}\otimes
\textrm{Ad}(\mc F_i),\theta_i)_{i=1}^r},
\end{equation*}
where $C=\max_{1\le i\le r} \theta_i^{-1}$, and  $q_n=\max \supp x_n$ for each $n$.  Find an infinite set $M$
of integers  and a sequence $(\mc G_i)_{i=1}^r$ of regular families such that for every $1\le i\le r$

\noindent (a) $\mr{Ad}(\mc F_i)\rest M$ is homogeneous on $M$,

\noindent (b)  $\iota (\mc G_i)=\iota(\mr{Ad}([M]^{\leq 2}\otimes \mr{Ad}(\mc F_i)\rest M))+1$, and

\noindent (c)  $ \mr{Ad}([M]^{\leq 2}\otimes \textrm{Ad}(\mc F_i))\rest M\con \mc G_i\rest M$

By Proposition \ref{xxsllki} we get
\begin{equation}
\nrm{\sum_{n\in M}a_n e_{q_n}}_{([\N]^{\leq 2}\otimes \textrm{Ad}(\mc F_i),\theta_i)_{i=1}^r} \le
\nrm{\sum_{n\in M}a_n e_{q_n}}_{(\mc G_i,\theta_i)_{i=1}^r}.
\end{equation}
By Theorem \ref{vmnjdgdjgjrre} there is some $N\con M$ and $C\ge 1$ such that
\begin{equation}
\nrm{\sum_{n\in N}a_n e_{q_n}}_{(\mc G_i,\theta_i)_{i=1}^r}\le C \nrm{\sum_{n\in N}a_n e_{q_n}}_{(\mc
G_{i_0},\theta_{i_0})},
\end{equation}
where $i_0$ is such that $(\iota (\mc G_{i_0}),\theta_{i_0})=\max_{<_\mathrm{T}}\conj{(\iota (\mc
G_i),\theta_i)}{1\le i\le r} $. Notice that $\iota(\mc G_{i_0})$ and $\iota(\mc F_{i_0})$ are both
infinite. By Corollary \ref{jfdfsfsdhsdrevbtghgg} we can find $P\con N$ such that
\begin{equation}
\nrm{\sum_{n\in P}a_n e_{q_n}}_{(\mc G_{i_0},\theta_{i_0})}\le 2\nrm{\sum_{n\in P}a_n e_{p_n}}_{(\mc
G_{i_0},\theta_{i_0})}
\end{equation}
where $p_n=\min \supp x_n$ for every $n\in P$. Since, by the choice of $M$, the $\mr{Ad}(\mc
F_{i_0})\rest M$ is homogeneous on $M$ we obtain that   by Proposition \ref{ceyhrytyyyy} that
\begin{equation}\label{gfwehgww}
\iota (\mc G_{i_0})\le 2 \iota([M]^{\le 2}\otimes\mr {Ad} (\mc F_{i_0})\rest M)+1 = 2
\iota([M]^{\le 2})\iota(\mr {Ad} (\mc F_{i_0})\rest M)+1 \le  8 \iota (\mc F_{i_0})+1
\end{equation}
Use now Proposition \ref{dfderjjwewww} to find an infinite subset $Q\con P$ and a regular family
$\mc B$ with the same index than $\mc F_{i_0}$ such that every $\mc B$-admissible sequence of
subsets of $\{p_n\}_{n\in Q}$ is $\mc F_{i_0}$-admissible.

Since   $\iota (\mc F_{i_0})$ is infinite and  $\mc G_{i_0}$ and $\mc B$ are regular, the
inequality (\ref{gfwehgww}) implies that
\begin{equation} \iota( [\N]^{\le 2} \otimes \mc G_{i_0})= 2\iota (\mc G_{i_0})\le 16 \iota (\mc F_{i_0})+2 < \iota (\mc F_{i_0}
)2= \iota (\mc B\otimes [\N]^{\le 2} ),
\end{equation}
so we can find an infinite $R\con Q$ such that
\begin{equation}
[\{p_n\}_{n\in R}]^{\le 2} \otimes \mc G_{i_0}\rest \{p_n\}_{n\in R} \con \mc B \otimes [\N]^{\le 2}.
\end{equation}
Hence, by Proposition \ref{xxsllki},
\begin{equation}
\nrm{\sum_{n\in R} a_{n}e_{p_n}}_{(\mc G_{i_0},\theta_{i_0})}\le 2 \nrm{\sum_{n\in R} a_{n}e_{p_n}}_{(\mc
B,\theta_{i_0})},
\end{equation}
while by Lemma \ref{b47}  we can find $S\con R$ such that
\begin{equation}
\nrm{\sum_{n\in S} a_n e_{p_n}}_{(\mc F_{i_0},\theta_{i_0})}\le 2 \nrm{\sum_{n\in S} a_n x_n}_{(\mc
F_{i_0},\theta_{i_0})}.
\end{equation}
Putting all these inequalities together we obtain
\begin{align}
\nrm{\sum_{n\in S}a_n e_{p_n}}_{(\mc G_{i_0},\theta_{i_0})} \le  & 2\nrm{\sum_{n\in S}a_n e_{p_n}}_{(\mc
B,\theta_{i_0})} \le  2\nrm{\sum_{n\in S}a_n e_{p_n}}_{(\mc
F_{i_0},\theta_{i_0})}\le \nonumber \\
& \le  4\nrm{\sum_{n\in S}a_n x_n}_{(\mc F_{i_0},\theta_{i_0})}\le 4\nrm{\sum_{n\in S}a_n x_n}_{(\mc
F_{i},\theta_{i})_{i=1}^r}.
\end{align}
So, $(x_n)_{n\in S}$ and $(e_{p_n})_{n\in S}\con T(\mc F_{i_0},\theta_{i_0})$ are equivalent, as desired.
\fprue
Recall from Definition \ref{fgjrioeeee} that  for a given compact and hereditary family $\mc F$ we set $
\ga(\mc F)$ and $n(\mc F)$ for $\iota(\mc F)$ and 1 respectively if $\mc F$ has finite index, and  for
$\om^{\om^{\ga}}$ and $n$ satisfying that $\om^{\om^{\ga}n}\le \al < \om^{\om^{\ga}(n+1)}$,  if $\mc F$ has
infinite index. Using this terminology we have the following

\begin{theorem*}\cite{ber-del},\cite{ber-pas}
Fix $(\mc{F}_i,\theta_i)_{i=1}^r$. Let $i_0$ be such that $(\mc F_{i_0},\theta_{i_0})=\max_{\le_\mr
T}\{(\mc F_i,\theta_i)\}_{i=1}^r$, and $\mc G$ be an arbitrary regular family  such that $\ga(\mc
G)=\ga(\mc F_{i_0})$. Then every normalized block sequence $(x_n)$ of $T[\fami]$ has a  block
subsequence $(y_n)_{n}$ equivalent to the basis of $T(\mc G,\theta_{i_0})$.
\end{theorem*}
Notice that the condition $\ga(\mc G)=\ga(\mc F_{i_0})$ above is equivalent to $\ga(\mc G)=\iota(\mc
G)=\iota(\mc F_{i_0})=\ga(\mc F_{i_0})$, and that $\mc G$ is equal to $[\N]^{\le \iota(\mc G)}$ in a tail
$\N/m$.  Now, in the same direction of this Theorem,

\cor
\label{kbmbnreyryy} Fix $(\mc{F}_i,\theta_i)_{i=1}^r$. Let $i_0$ be such that $(\mc
F_{i_0},\theta_{i_0})=\max_{\le_\mr T}\{(\mc F_i,\theta_i)\}_{i=1}^r$. Suppose that $\mc G$ is an
arbitrary compact and hereditary family. If  $\ga(\mc G)=\ga(\mc F_{i_0})$, then every normalized
block sequence $(x_n)$ of $T[\fami]$ has a   subsequence $(x_n)_{n\in M}$ equivalent to the
subsequence $(e_{\min \supp x_n})_{n\in M}$ of the basis of $T(\mc G,\theta_{i_0}^{n(\mc B)/n(\mc
F_{i_0})})$.
\fcor
\prue
By Theorem \ref{infiniteindex} it is enough to have the conclusion for subsequences of the basis of
$T(\mc F_{i_0},\theta_{i_0})$, and by Proposition \ref{dfewekmmvv2} we may assume that $\mc
F_{i_0}$ and $\mc G$ are both regular families. Set $\iota(\mc F_{i_0})=\om^{\om^\al m+\be} n
+\de$,  $\iota(\mc G)=\om^{\om^{\bar{\al}} \bar{m}+\bar{\be}} \bar{n} +\bar{\de}$ be  canonical
decompositions This is possible since $\ga(\mc F_{i_0})=\ga(\mc B)$ is infinite. Moreover
$\bar{\al}=\al$. Using
\begin{equation*}
\om^{\om^\al m}\le \iota(\mc F_{i_0})=\om^{\om^\al m+\be} n +\de< \om^{\om^\al m+\be+1},
\end{equation*}
and the corresponding inequality for $\mc G$, by Theorem \ref{surprise} we may assume that
$\iota(\mc F_{i_0})=\om^{\om^\al m}$, and $\iota(\mc G)=\om^{\om^\al \bar{m}}$ Now the result
follows from the application of Proposition \ref{powerisnot}  to the families $\mc F_{i_0}$ and
$\mc G$.
\fprue

In particular for Schreier families we obtain
\cor
\label{thf} Fix $(\mc{F}_i,\theta_i)_{i=1}^r$ such that at least one of the families has infinite
index. Let $i_0$ be such that $(\mc F_{i_0},\theta_{i_0})=\max_{\le_\mr T}\{(\mc
F_i,\theta_i)\}_{i=1}^r$, and set $\iota(\mc F_{i_0})=\om^{\om^\al k+\de}m+\ga$ in canonical form.
Then every normalized block sequence $(x_n)$ of $T[\fami]$ has a subsequence $(x_n)_{n\in M}$
equivalent to the subsequence $(e_{\min \supp x_n})_{n\in M}$ of the basis of $T(\mc
S_{\om^{\al}},\theta_{i_0}^{1/k})$.\qed
\fcor

The last result of the section  concerns  equivalence  of block sequences in the spaces $T(\mc
S_{\xi}, \theta)$.

\prop
Let $(x_n)$, $(y_n)$ be two normalized block sequences in the space $T(\mc S_{\xi},\theta)$ be such
that $x_n<y_n <x_{n+1}$ ($n\in \N$). Then $(x_n)$ and $(y_n)$ are $24\theta^{-2}$-equivalent.
\fprop
\prue
For the proof  we shall use the following two relations concerning  the Schreier families $\mc
S_{\xi}$, and infinite subsets $N$ of integers with $\min N\geq 3$.
\begin{align}
[N]^{\leq 3}\otimes \mc S_{\xi}\subseteq & \mc S_{\xi}\otimes [N]^{\leq 2}\label{ee1}\\
[N]^{\leq 3}\otimes (\mc S_{\xi}\oplus  [N]^{\leq 1})\subseteq & \mc S_{\xi}\otimes [N]^{\leq
3}\label{ee2}
\end{align}
The proof of these two relations follows easily by induction on $\xi$. We show now that  a normalized block
sequence $(x_n)$  is equivalent to the subsequence $(e_{p_n})_{n}$ of the basis, $p_n=\min\supp x_n$,  and
this implies the result. Without loss of generality we may assume that $p_n\geq 3$ for every $n$. It follows
easily form the spreading property of the families $\mc S_{\xi}$ that
$$
\nrm{\sum_n a_n e_{p_n}}_{(\mc S_{\xi},\theta)}\leq \nrm{\sum_n a_n x_{n}}_{(\mc S_{\xi},\theta)}.
$$
For the reverse inequality,  by Proposition \ref{dwbvsaadfsdfsee} we get
\begin{equation*}
\nrm{\sum_n a_n x_n}_{(\mc S_{\xi},\theta)}\le \theta^{-1} \nrm{\sum_n a_n e_{q_n}}_{([\N]^{\leq
2}\otimes \mc S_{\xi},\theta)},
\end{equation*}
where   $q_n=\max \supp x_n$ for each $n$. By \eqref{ee1}  and Proposition \ref{xxsllki} we get
$$
\nrm{\sum_n a_n e_{q_n}}_{([\N]^{\leq 2}\otimes \mc S_{\xi},\theta)} \leq 2\nrm{\sum_n a_n
e_{q_n}}_{(\mc S_{\xi},\theta)}.
$$
As in the proof of Corollary \ref{jfdfsfsdhsdrevbtghgg} we get that
$$
\nrm{\sum_n a_n e_{q_n}}_{(\mc S_{\xi},\theta)}\leq \nrm{\sum_n a_n e_{p_n}}_{(\mc S_{\xi}\oplus
[\N]^{\leq 1},\theta)}.
$$
Now by \eqref{ee2} and again Proposition \ref{xxsllki} we get that
$$
\nrm{\sum_n a_n e_{p_n}}_{(\mc S_{\xi}\oplus [\N]^{\leq 1},\theta)}\leq 3\nrm{\sum_n a_n
e_{p_n}}_{(\mc S_{\xi},\theta)}.
$$
and this completes the proof.
\fprue

\subsection{Incomparability} The goal here is to turn the implication presented in Corollary
\ref{kbmbnreyryy} into an equivalence.  So we are now going to  deal  with the incomparability of
the Tsirelson-type spaces. The main tool to distinguish two such spaces are the special convex
combinations, introduced in \cite{arg-del}. The following lemma provides the existence of the
special convex combinations,  in a more general setting than the one in \cite{arg-del}, and it is a
version of the well known Pt\'{a}k's Lemma (see \cite{arg-tod} for a proof).
  \lema \label{ptak}
Suppose that $\mc F_0$ and $\mc F_1$ are two regular families  with indexes $\iota (\mc
F_i)=\om^{\al_i} n_i +\be_i$, $\al_i>0$, $n_i\in \N$ $\be_i<\om^{\al_i}$ ($i=0,1$). If
$\al_0<\al_1$, then for every $\vep>0$ there is a convex mean $\mu$ such that $\supp \mu\in \mc
F_1$  and such that $\sup_{t\in \mc F_0}\sum_{n\in t}\mu(n)<\vep$.
 \flema
The first case where the spaces are going to be totally incomparable
is if the index of one of the families is at least the $\om$-power
of the other.
 \lema\label{knryyyyyuuii}
Suppose that $\mc F_0,\mc F_1$  are two regular families such that $\iota (\mc F_0)^\om \le \iota
(\mc F_1)$. Then $ T(\mc F_0,\theta_0) $ and $ T(\mc F_1,\theta_1) $ are totally incomparable. \qed
 \flema
 \prue
Suppose that the desired result does not hold.
 By standard arguments we may assume that
there exists a normalized block sequence $(x_{n})_{n}\in T(\mc F_i,
\theta_i)$ equivalent to a normalized block sequence $(z_{n})_{n}$
of $T(\mc F_j, \theta_j)$, $j\ne i$. By Theorem \ref{infiniteindex}
passing to subsequences if necessary we may assume that
$(x_{n})_{n}$ is equivalent to a subsequence $(e_n)_{n\in M_i}$ of
the natural basis $(e_n)$ of $ T(\mc F_i,\theta_i) $ and that
$(z_n)$ is equivalent to a subsequence $(e_n)_{n\in M_j}$ of the
natural basis $(e_n)$ of $ T(\mc F_j,\theta_j) $.

 For $k=0,1$, let $\vphi_k:M_k\to \N$ be the unique
order-preserving onto mapping between $M_k$ and $\N$. Note that for $k=0,1$ the family
$\phi_k^{-1}\mc F_k$ is regular on $M_{k}$, $\iota(\vphi_k^{-1}\mc F_k)=\iota (\mc F_k)$ and
$(e_n)_{n\in M_k}\con T(\mc F_k,\theta_k)$ is 1-equivalent to $(e_n)_{n\in \N}\con
T(\vphi_k^{-1}\mc F_k,\theta_k)$. So, without loss of generality, we may assume that $M_1=M_2=\N$.
So, we are supposing that $(e_n)\con T(\mc F_0,\theta_0) $ is, say, $C$-equivalent to $(e_n)\con
T(\mc F_1,\theta_1) $ i.e. for every scalars $(a_n)$,
\begin{equation}\label{jfjhfghweeemm0}
\frac1C\nrm{\sum_n a_n e_n}_{(\mc F_0,\theta_0)}
 \le
\nrm{\sum_n a_n e_n}_{(\mc F_1,\theta_1)}
 \le
C\nrm{\sum_n a_n e_n}_{(\mc F_0,\theta_0)}.
\end{equation}
Let $l\in \N$ be such that $\theta_0^{l}<\theta_1/(2C)$. By our hypothesis over the indexes, $\iota (\mc
F_0)^l<\iota (\mc F_1)$. So, by  Lemma \ref{ptak} there is some convex mean $\mu$ such that
\begin{equation*}
\text{$\supp \mu \in \mc F_1$, and   $\sum_{n\in t}\mu(n)< \frac{\theta_1}{2C}$ for every $t\in \mc
F_0^{\otimes(l-1)}$}.
\end{equation*}
Observe that every $\phi\in K(\mc F_0,\theta_0)$ has a decomposition $\phi=\phi_0+\phi_1$, where $\supp
\phi_0\in  \mc F_0^{\otimes(l-1)}$,  $\nrm{\phi_1}_{\infty}\le \theta^{l}$ and $\supp \phi_0\cap \supp
\phi_1=\buit.$ So, for every $\phi\in K(\mc F_0,\theta_0)$,
\begin{equation*}
\begin{split}
|\phi(\sum_{n\in s}\mu(n) e_n)|
& =|\phi_0(\sum_{n\in s}\mu(n) e_n)+
\phi_1(\sum_{n\in s}\mu(n) e_n)|
\\
 &\le
\sum_{n\in \supp \phi_0\cap s}\mu(n) + \nrm{\phi_1}_\infty
\sum_{n\in s}\mu(n)< \frac{\theta_1}{2C}+
\theta_0^l<\frac{\theta_1}{C},
\end{split}
\end{equation*}
while
\begin{equation*}
\nrm{\sum_{n\in s}\mu(n)e_n}_{ (\mc F_1,\theta_1) }\ge
\theta_1\sum_{n\in s}\mu(n)=\theta_1,
\end{equation*}
and so, by (\ref{jfjhfghweeemm0}),
\begin{equation*}
\theta_1\le \nrm{\sum_{n\in s}\mu(n)e_n}_{  (\mc F_1,\theta_1) } \le
C \nrm{\sum_{n\in s}\mu(n)e_n}_{  (\mc F_0,\theta_0)
}<C\frac{\theta_1}C,
\end{equation*}
a contradiction.
   \fprue
The second case of totally incomparability we consider is when the two families have the same
index, but the corresponding $\theta$'s are different.
\lema\label{sameindex} Suppose that $\mc F_0$ and $\mc F_1$ are two
regular families with same index, and suppose that $\theta_0\neq \theta_1$, and
$\max\{\theta_0,\theta_1\}>1/\iota (\mc F_0)$, where by convention, $1/\al=0$ if $\al$ is an
infinite ordinal. Then the corresponding spaces $ T(\mc F_0,\theta_0) $ and $ T(\mc F_1,\theta_1) $
are totally incomparable.
\flema
 \prue
Suppose first than $\iota(\mc F_0)=\iota(\mc F_1)$ is finite. Then $T(\mc F_0,\theta_0)$ and $T(\mc
F_1)$ are isomorphic to different classical spaces $c_0$ or $\ell_p$ ($p > 1$), and the conclusion
of the Lemma trivially holds.

 Suppose that $\iota (\mc F_0)=\iota (\mc F_1)$ is infinite. As in previous lemma, we may assume
that $\mc F_0=\mc F_1=\mc F$ and that $(e_n)\con T(\mc F,\theta_0)$ is $C$-equivalent to $(e_n)\con
T(\mc F,\theta_1)$, i.e. for every scalars $(a_n)$,
\begin{equation}\label{jfjhfghweeemm}
\frac1C\nrm{\sum_n a_n e_n}_{ (\mc F,\theta_0)}
 \le
\nrm{\sum_n a_n e_n}_{ (\mc F,\theta_1)}\le C\nrm{\sum_n a_n e_n}_{
(\mc F,\theta_0)}.
\end{equation}
Suppose that $\theta_0<\theta_1$.  Let $l\in \N$, $l>1$ be such that $(\theta_1/\theta_0)^{l}>2C$. Let
$(a_n)_{n\in s}$ be a convex mean such that $s\in \mc F^{\otimes(l)}$ and   $\sum_{n\in
t}a_n<\theta_1^l/(2C)$  for every $t\in \mc F^{\otimes(l-1)}$. As before, any functional $\phi\in K(\mc
F,\theta_0)$ is decomposed $\phi=\phi_0+\phi_1$, $\supp \phi_0\cap \supp \phi_1=\buit$, $\supp \phi_0\in \mc
F^{\otimes(l-1)}$ and $\nrm{\phi_1}_\infty\le \theta_0^{l}$. Then
\begin{align}
|\phi(\sum_{n\in s}a_n e_n)|
 =|\phi_0(\sum_{n\in s}a_n e_n)+
\phi_1(\sum_{n\in s}a_n e_n)|
 < \frac{\theta_1^l}{2C}+ \theta_0^l<
\frac{\theta_1^l}{C}
\end{align}
Finally, by (\ref{jfjhfghweeemm}),
\begin{align}
\theta_1^l\le \nrm{\sum_n a_n e_n}_{(\mc F,\theta_1)}\le
C\nrm{\sum_n a_n e_n}_{(\mc F,\theta_0)}< C\frac{\theta_1^l}{C},
\end{align}
a contradiction.
  \fprue
  \subsection{Main result}
We collect  in a single result the facts we have got so far.

\teor[Classification theorem]\label{maintheorem}
Fix two sequences $(\mc F_i,\theta_i)_{i=1}^r$ and  $(  {\mc G}_i,{\eta}_i)_{i=1}^{s}$ of pairs of
compact and hereditary families and real numbers in $(0,1)$. Let $1\le i_0\le r$ and $1\le {j}_0\le
s$ be such that $(\mc F_{i_0},\theta_{i_0})=\max_{\le_\mr T} \conj{(\mc F_i,\theta_i)}{1\le i\le
r}$, and $({\mc G}_{j_0},{\eta}_{j_0})=\max_{\le_\mr T} \conj{({\mc G}_i, {\eta}_i)}{1\le i\le s}$.
The following are equivalent:

\noindent (a)   Either

\noindent \, (a.1) $\ga(\mc F_{i_0}),\ga(\mc G_{j_0})\ge \omega$, $\ga(\mc F_{i_0})=\ga(\mc
G_{j_0})$ and $\theta_{i_0}^{n(\mc G_{j_0})}=\eta_{j_0}^{n(\mc F_{i_0})}$, or else

\noindent \, (a.2) both $\mc F_{i_0},\mc G_{j_0}$ have finite index, and  either

\noindent \quad (a.2.1) $\theta_{i_0}\ga(\mc F_{i_0}), \eta_{i_0}\ga(\mc G_{j_0})\le 1$, or else

\noindent \quad (a.2.2) $\log_{\ga(\mc F_{i_0})}\theta_{i_0}=\log_{\ga(\mc G_{j_0})}\eta_{j_0}$.

\noindent (b) Every closed infinite dimensional subspace of $ T[(\mc F_i,\theta_i)_{i=1}^r] $
contains a subspace isomorphic to a subspace of $T[(\mc G_i,\eta_i)_{i=1}^s]$.

\noindent (c) For every regular family $\mc B$ such $\ga(\mc B)=\ga(\mc G_{j_0})$ and every
normalized block sequence of $ T[(\mc F_i,\theta_i)_{i=1}^r] $ there is a block subsequence
(subsequence if $\mc G_{i_0}$ has infinite index) equivalent to a subsequence of the natural basis
of $ T(\mc B,\eta_{j_0}^{n(\mc B)/n(\mc G_{j_0})}) $.

\fteor
\prue
(b) implies (c). Fix a regular family $\mc B$  with same index than   $\mc G_{j_0}$, and fix a
normalized block sequence $(x_n)$ of $ T[\fami]$. By (b), there is some block sequence $(y_n)$ of
$(x_n)$ which is equivalent to a semi normalized block sequence $(z_n)$ of $ T[(\mc
G_i,\eta_i)_{i=1}^s] $. By  Corollary \ref{kbmbnreyryy}, we can find a further block subsequence
$(w_n)$ of $(z_n)$ which is equivalent to a subsequence of the natural basis of $T(\mc
B,\eta_{j_0}^{1/n(\mc G_{j_0})})$, as desired.

\medskip

(c) implies (a).    First of all, fix a regular family $\mc C$ with index $\ga(\mc F_{i_0})$. By
Corollary \ref{kbmbnreyryy} we know that $T(\fami)$ is saturated by subsequences of the basis of
$T(\mc C,\theta_{i_0}^{1/n(\mc F_{i_0})})$. Notice that (c) implies that  $T(\mc
B,\eta_{j_0}^{1/n(\mc G_{j_0})})$ and $T(\mc  C,\theta_{i_0}^{1/n(\mc F_{i_0})})$ are not totally
incomparable. Suppose first than $\mc G_{j_0}$ has finite index. Lemma \ref{knryyyyyuuii} gives
that $\mc F_{i_0}$ has
 also finite index, and in particular $n(\mc F_{i_0})=1$. Now (a.2) follows from the properties of $\ell_p$'s and
 $c_0$.

Assume now that $\mc G_{j_0}$ is infinite. In this case Lemma   \ref{knryyyyyuuii}  implies that
$\ga(\mc F_{i_0})=\ga(\mc G_{j_0})$. It follows,  by Corollary \ref{kbmbnreyryy} that $T[\fami]$ is
saturated by subsequences of $T(\mc B,\theta_{i_0}^{1/n(\mc F_{i_0})})$. Hence   $T(\mc
B,\theta_{i_0}^{1/n(\mc F_{i_0})})$  and $T(\mc B,\eta_{j_0}^{1/n(\mc G_{j_0})})$ are not totally
incomparable, so by Lemma \ref{sameindex}, $\theta_{i_0}^{1/n(\mc F_{i_0})}=\eta_{i_0}^{1/n(\mc
G_{j_0})}$.

(a) implies (b) follows from Corollary \ref{kbmbnreyryy}.
\fprue

\nota \mbox{}

\noindent 1. If the families $\mc F$ are compact but not necessarily hereditary, Theorem
\ref{maintheorem} is also true. The main observation is that  if $\mc F$ is arbitrary compact
family, there is some infinite set $M$ such that $\mc F[M]=\conj{s\cap M}{s\in \mc F}$ is
hereditary (see \cite{arg-hand}). This fact when applied to the family $\mathrm{Ad}(\mc F)$ of $\mc
F$-admissible sets guarantees to follow the arguments we use for the case of hereditary families,
starting with Proposition \ref{dwbvsaadfsdfsee}.

\noindent  2.  The problem of classification of full mixed Tsirelson spaces $T[(\mc
F_i,\theta_i)_{i=0}^\infty]$ seem rather unclear. There are several obstacles if someone wants to
extend the techniques presented in this paper to the general case.
\fnota

\end{document}